\documentclass[10pt]{beamer}

\usepackage{multimedia}
\usepackage{relsize}
\usepackage[utf8x]{inputenc}
\usepackage{subfigure}

\usepackage[english]{babel}
\usepackage{amsmath,amssymb}
\usepackage[dvips]{epsfig}
\usepackage{pstricks,pst-node}
\newtheorem{prop}[]{Proposition}

\title[Robustification of Opimal Controls.]
{On Optimized Feedback Control and the Robustification of Opimal Controls
}

\author[Martin Gugat]
{Martin Gugat} \institute[FAU]
{\inst{}FAU: Friedrich-Alexander-Universit\"at Erlangen-N\"urnberg}

\date[]{
{\bf CDPS 2013}:
8th Workshop on Control of Distributed Parameter Systems,
July 1- 5, 2013, University of Craiova
}

\begin{document}

\maketitle






\section[Outline]{}
\frame{\tableofcontents}


\begin{frame}[t]

\titlepage
\begin{figure}[!h]

\includegraphics[width=5cm]{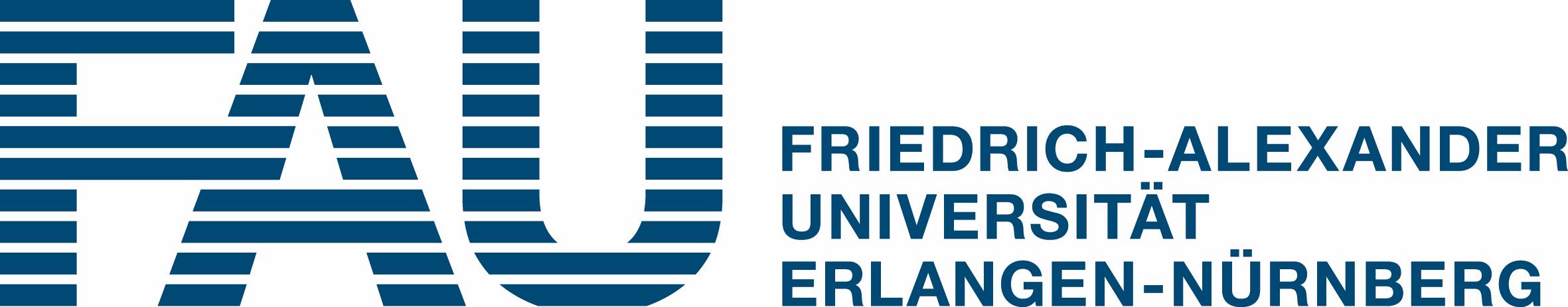}


\end{figure}

\end{frame}

%
%
%
%
%
%
%
%
%
%

%
%
%
%
%
%
%
%
%
%
%
%
%
%
%
%
%
%
%
%
%



\begin{frame}

{\Huge \bf Optimal Boundary Control of the Wave Equation}

\end{frame}


\section{Optimal Dirichlet Boundary Control}


\begin{frame}

{\Huge {\bf Optimal Dirichlet Boundary Control}

$$ \fbox{$y(t,1)= u (t)$}$$

}

\end{frame}


\begin{frame}\frametitle{The Problem of Optimal Exact Control: The 1d-case}

\begin{itemize}

\item<1->
Let the final time  $T=2k$ with a  natural number $k$ be given.

\item<2->
We consider the {\em wave equation}
on $[0,T]\times[0,1]$.

\item<3->
Initial position $y_0 \in L^2(0,1)$.

\item<4->
Initial velocity $y_1 \in H^{-1}(0,1)$.

\item<5->

%
$$
\label{EC} {\bf (EC)}
\left\{
\begin{array}
{l}
{\rm minimize} \;
 \|u\|_{L^2(0,T)}^2
\; {\rm subject } \;{\rm to}\;
\\
\\
y(0,x)= y_0(x),\; y_t(0,x)= y_1(x),\; x\in (0,1)
\\
\\
y(t,0)= 0,\;  \fbox{$ y(t,1)= u(t),$}\; t\in (0,T)
\\
\\
y_{tt}(t,x)= y_{xx}(t,x),\;(t,x)\in (0,T) \times (0,1)
\\
\\
y(T,x)= 0,\; y_t(T,x)= 0,\; x\in (0,1).
\end{array}
\right.
$$
%
\end{itemize}

\end{frame}



\begin{frame}\frametitle{Solution of  Problem (EC)}
\begin{itemize}

\item<1->
%
%
%
Problem ${\bf EC}$ has a solution $u$ that
is uniquely determined.

\item<2->
The optimal control $u_\ast$  is $2$ periodic.

\item<3->
$$
u_\ast(t) =
\left\{
\begin{array}{ll}
\frac{1}{T} \left( -  \int_0^{1-t} y_1(s)\, ds + r  + y_0(1-t)\right), &  t\in (0,1) \\
\frac{1}{T} \left( - \int_0^{t-1} y_1(s)\, ds + r  -y_0(t-1)\right), &  t\in (1,2)
\end{array}
\right.
$$
with $r= \int_0^1 \int_0^t y_1(s)\,ds \, dt.$
\\
{\tiny \em M. Gugat, G. Leugering, G. Sklyar: Lp optimal boundary control for the wave equation,
SICON 2005}

%
%




\end{itemize}

\end{frame}


\begin{frame}\frametitle{Penalization of the  exact end conditions}

\begin{itemize}

\item<1->
The exact end conditions in {\bf EC} can be replaced  by a
\\
{\bf    nondifferentiable  penalty term} in the objective function:
$$\min \;
\frac{1}{\gamma}\; \|u\|_{L^2(0,T)}^2 +
\sqrt{
 \|y(T,\cdot)\|_{L^2(0,1)}^2
+ \|Y\|_{L^2(0,1)}^2},
$$
$$
 Y(0)   = -  \int_0^1 \int_0^x y_t(t,\,z)\, dz\,dx  ,\;   Y'(x)=  y_t(T,x).
$$
\item<2->

For $\gamma \geq \frac{2}{\sqrt{k}} \|u_\ast\|_{L^2(0,T)}$,
this problem also has the solution $u_\ast$ of {\bf EC}.


\vspace*{5mm}

%
{\em M. Gugat: Penalty Techniques for State Constrained Optimal Control Problems with the Wave Equation,
SICON 2009}


\item<3->
This problem has a solution also for small $T$.

\end{itemize}

\end{frame}

\begin{frame}\frametitle{Example}
\begin{itemize}

\item<1->
Let $y_0(x)=x$, $y_1(x)=0$.

\item<2->
We get the optimal control
%
$ u_\ast(t) = \frac{1}{T} \left( 1-t\right),  t\in (0,2) $.

\hspace*{1cm}
\includegraphics[width=8cm,height=3cm]{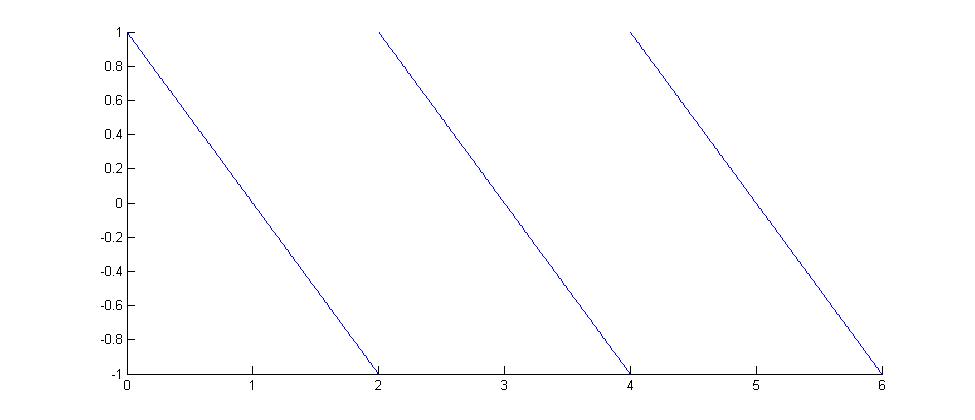}

\item<3->
Thus  if $T>2$, we have a {\bf jump} at time $t=2$!

\item<4->
Hence also for continuous data, the optimal state
for {\sc Dirichlet} control is in general
{\em discontinous}.
{\bf Continuity} is an {\bf additional constraint},
\\
{\tiny see  \sc M. Gugat; \em  Optimal boundary control of a string
to rest in finite time with continuous state, \sf ZAMM,  86 (2006)
pp. 134-150. }

\item<6-> To do this, we need $y_0\in H^1(0,1)$, $y_1\in L^2(0,1)$.

\end{itemize}

\end{frame}


\begin{frame}\frametitle{Continuous states}

\begin{itemize}

\item<1->
The following optimal control problem admits only continuous states:
$$
 { \cal P}
\left\{
\begin{array}{l}
 {\rm minimize} \; \|(u_0',\,u_1')\|_{2,(0,T)}
\; {\rm subject } \;{\rm to}\;
\\
\\
u_0, u_1 \in  H^1[0,T]
\\
\\
 y(0,x)= y_0(x),\; y_t(0,x)= y_1(x),\; x\in (0,1)
 \\
 \\
\fbox{ $y(t,0)= u_0(t),\; y(t,1)= u_1(t),\; t\in [0,T]$}
\\
\\
y_{tt}(t,x)=  y_{xx}(t,x),\;(t, x)\in (0,T) \times (0,1)
\\
\\
y(T,x)= 0,\; y_t(T,x)= 0,\; x\in (0,1)
\\
\\
y_0(0) = u_0(0), \, y_0(1) = u_1(0), \, 0  = u_0(T),\,
 0  = u_1(T).
 \end{array}
 \right.
$$
In the last line you see $C^0$--compatibility conditions.

\end{itemize}

\end{frame}


%
%
%
%
%
%
%
%
%
%
%
%
%
%
%
%
%
%
%
%
%
%
%
%
%
%
%
%
%
%
%
%
\begin{frame}\frametitle{Continuous states}
Let $T=2$, $y_0(x)=-1$ and $y_1(x) = 0$.
\\ Optimal controls:
$u_0(t) = u_1(t) = -1 + t/2$.
%


\includegraphics[width=0.75\textwidth]{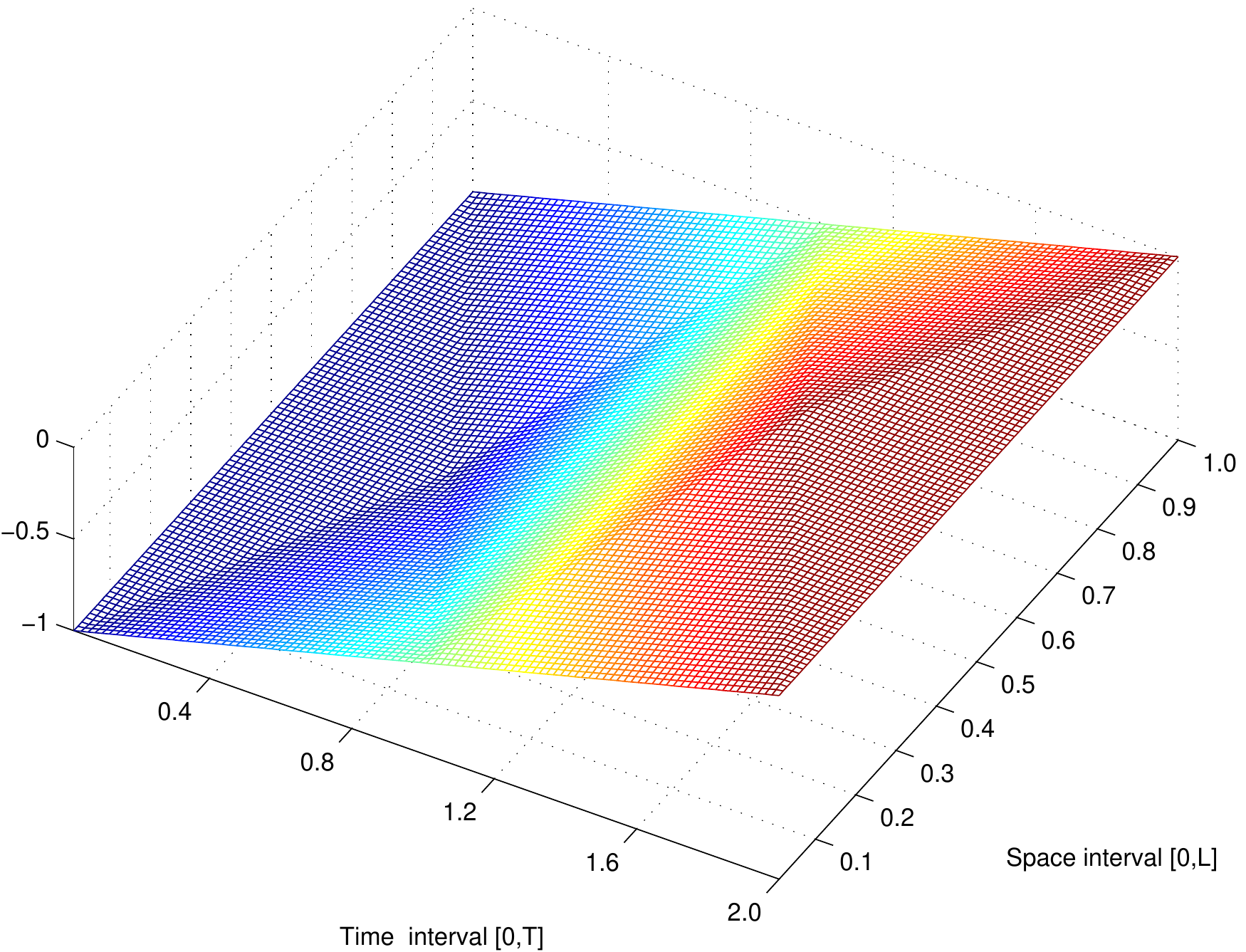}

\end{frame}


\begin{frame}\frametitle{Continuous states}
\begin{itemize}

\item<1-> With  {\sc Neumann} control, {\bf Continuity} is {\bf not}
an {\bf additional constraint}!

\item<2-> We will come to this later!

Let us first look at the $L^\infty$-case:

Do we get bang-bang controls?

\end{itemize}

\end{frame}




\frame{ \frametitle{$L^\infty$-case: Weakness of the bang-bang principle}
\begin{itemize}

\item<1-> $y_0\in L^\infty(0,1)$, $y_1 \in W^{-1,\infty}(0,1)$.
$$
 {\bf (DEC\infty)}
 \left\{
\begin{array}{l}
\min \; \frac{1}{2} \|u\|_{L^\infty(0,T)}^2
\; {\rm subject } \;{\rm to}\;
\\
y(0,x)= \sin(x \pi) ,\; y_t(0,x)= 0,\; x\in (0,1)
\\
y(t,0)= 0,\; y(t,1)= u(t),\; t\in (0,T)
\\
y_{tt}(t,x)= y_{xx}(t,x),\;(t,x)\in (0,T) \times (0,1).
\\
y(T,x)= 0,\; y_t(T,x)= 0,\; x\in (0,1).
\end{array}
\right.
$$

\item<2->
For $T=2$ an optimal control is
$$u(t)= \frac{1}{2} \sin(t \pi).$$
All admissible controls have the form $u(t)+const$, so
there is no admissible bang-bang control.

\item<3->
Let $T= 2k$.
States that can be reached by bang-bang-off controls:
$$y(x,T) \in y_0(x) + \|u\|_{\infty, (0,T)} \{-2k,\, -2k + 1,..., 2k -1,...,2k\}.$$


{\em M. Gugat, G. Leugering:
$L^\infty$ Norm Minimal Control of the wave equation: On the weakness of the bang--bang principle,
ESAIM: COCV 14
(2008)}

\end{itemize}

} 

\begin{frame}

{
\Huge
\bf
Now:
{\bf Neumann boundary control}

$$ \fbox{$y_x(t,1)= u (t)$}$$

}

\end{frame}




\section{Optimal Neumann Boundary Control}



\frame{
\frametitle{The Problem of optimal exact control: Neumann}

\begin{itemize}
\item<1->
Let $y_0 \in H^1(0,1)$, $y_1 \in L^2(0,1).$
\item<2->

$$
{\bf (EC)}
\left\{
\begin{array}
{l}
{\rm minimize} \;
 \|u\|_{L^2(0,T)}^2
\; {\rm subject } \;{\rm to}\;
\\
\\
y(0,x)= y_0(x),\; y_t(0,x)= y_1(x),\; x\in (0,1)
\\
\\
y(t,0)= 0,\; \fbox{$y_x(t,1)= u(t)$,}\; t\in (0,T)
\\
\\
y_{tt}(t,x)= y_{xx}(t,x),\;(t,x)\in (0,T) \times (0,1)
\\
\\
y(T,x)= 0,\; y_t(T,x)= 0,\; x\in (0,1).
\end{array}
\right.
$$

\end{itemize}

} 






\frame{ \frametitle{The method of characteristics: The key to the problem }



D'Alembert:
{\em Recherches sur la courbe que forme une corde tendue mise en vibration,
Mem. Acad. Sci. Berlin 3, 214-219, (1747).}

\vspace*{-3cm}

\includegraphics[width=9cm,angle=270]{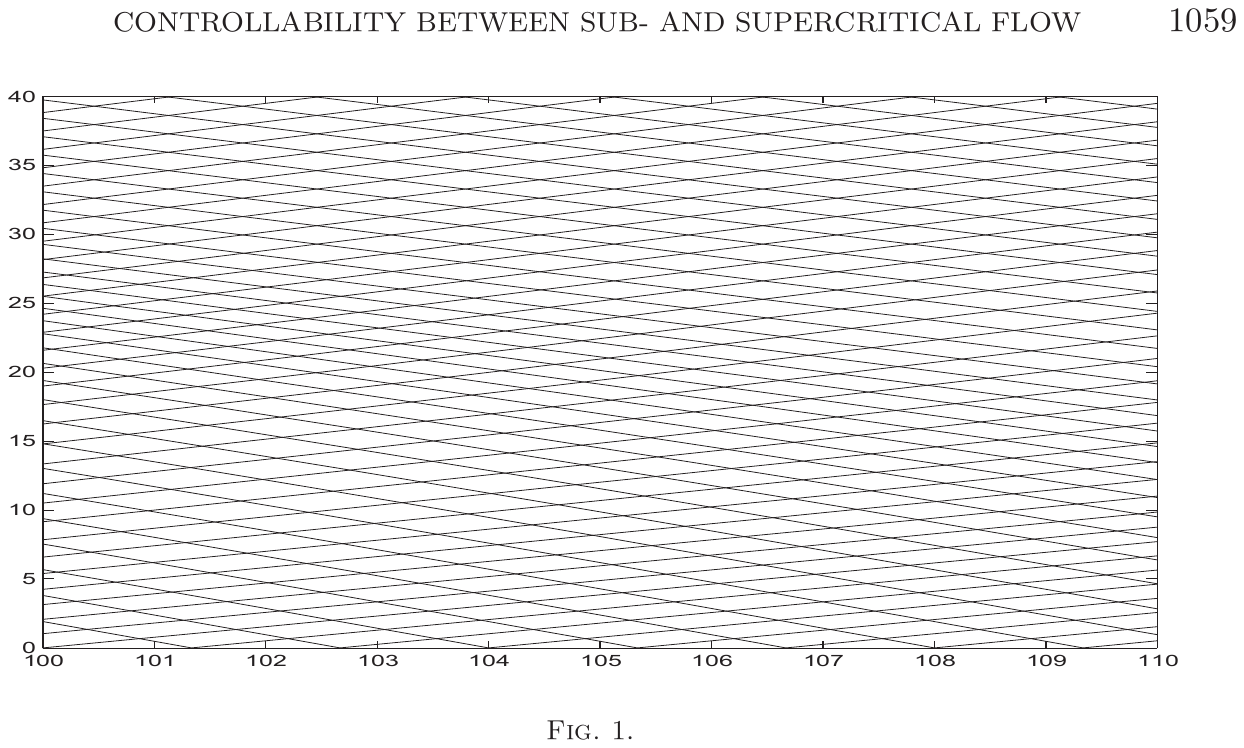}

\vspace*{-3cm}

D'Alembert's solution has the form
$$y(t,x)= \alpha(x+t) + \beta(x-t).$$


From the initial conditions for $t\in (0,1)$:

$\alpha(t)= \frac{1}{2} \left(y_0(t) + \int_0^t y_1(s)\, ds \right)+C,$
$\beta(t)= \frac{1}{2} \left(y_0(t) - \int_0^t y_1(s)\, ds \right)-C.$

}

\frame{ \frametitle{The optimal Neumann control}

\begin{itemize}

\item<1->



{\bf Theorem}
[{\small \em Gugat 2013}]
Let $T=K+1$ be even.

\item<2->
Then the optimal control is 4--periodic, with
$$
u(t) =
\left\{
\begin{array}{ll}
\frac{2}{T} \beta'(1-t)  =  \frac{1}{T}  \left(y_0'(1-t) - y_1(1-t)\right), &  t\in (0,1) \\
\frac{2}{T} \alpha'(t-1) =  \frac{1}{T}  \left(y_0'(t-1) + y_1(t-1)\right), & t\in (1,2).
\end{array}
\right.
$$




\item<3->

For $k\in \{0,1,...,(K-1)/2\}$, $t\in (0,2)$ we have:
$$u(t + 2k)= (-1)^k u(t).$$


\item<4->
\fbox{ Moving horizon idea:} At each moment, it is  best to use $u_\ast(0)$ 
\\
with the current state as initial data.

\item<5->
With the \fbox{ moving horizon idea} we get the feedback law:
%

$$
y_x(t,1)= \frac{1}{T-1}
 \left(
 - y_t(t,1)
\right)
$$

\item<6->
This is a well-known exponentially  stabilizing feedback!

\end{itemize}

} 

\subsection{Example: Solution of {\bf (EC)} }


\frame{ \frametitle{Example: Optimal Neumann Control}

Let $y_0(x) = 4\sin(\frac{\pi}{2} x)$, $y_1(x)=0$.
\vspace*{5mm}
Then $\alpha(x) = \beta(x) = 2\sin(\frac{\pi}{2} x)$.
\vspace*{5mm}

We obtain the optimal control
$$
u(t) =
\left\{
\begin{array}{ll}
\frac{2}{T}  \pi \cos(\frac{\pi}{2} (1 - t)), &  t\in (0,1);
\\
\\
\frac{2}{T}  \pi \cos(\frac{\pi}{2} (t-1)) , & t\in (1,2).
\end{array}
\right.
$$

By continuation we get
$$
\fbox{
$
u(t) = \frac{2}{T}  \,\pi \,\cos\left(\frac{\pi}{2} (t-1)\right).
$
}
$$

} 

\frame{ \frametitle{
Example: Minimal Control Time $T=2$:}

Optimal state for the minimal control time $T=2$:
%
\includegraphics[width=5cm]{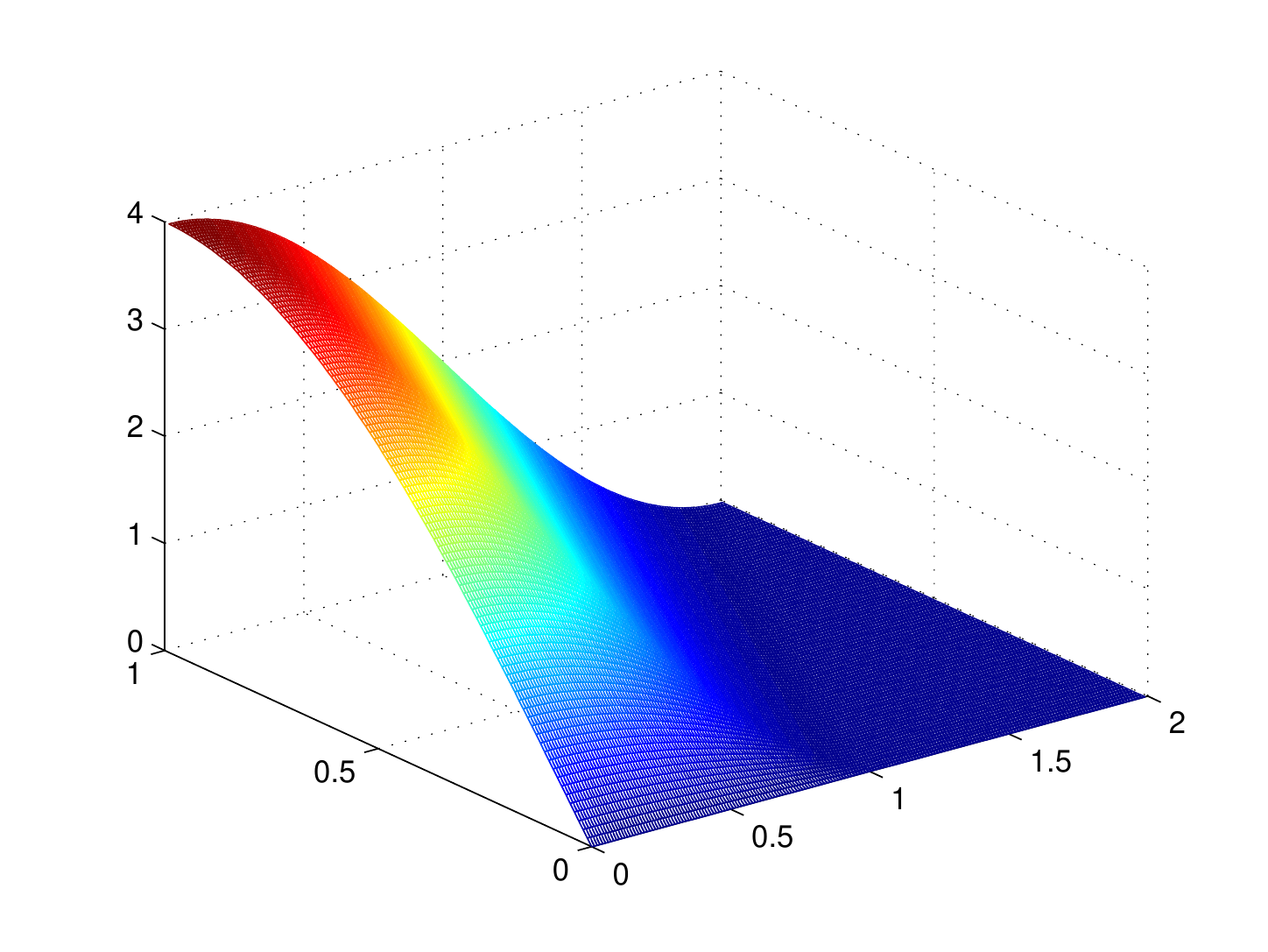}
\includegraphics[width=5cm]{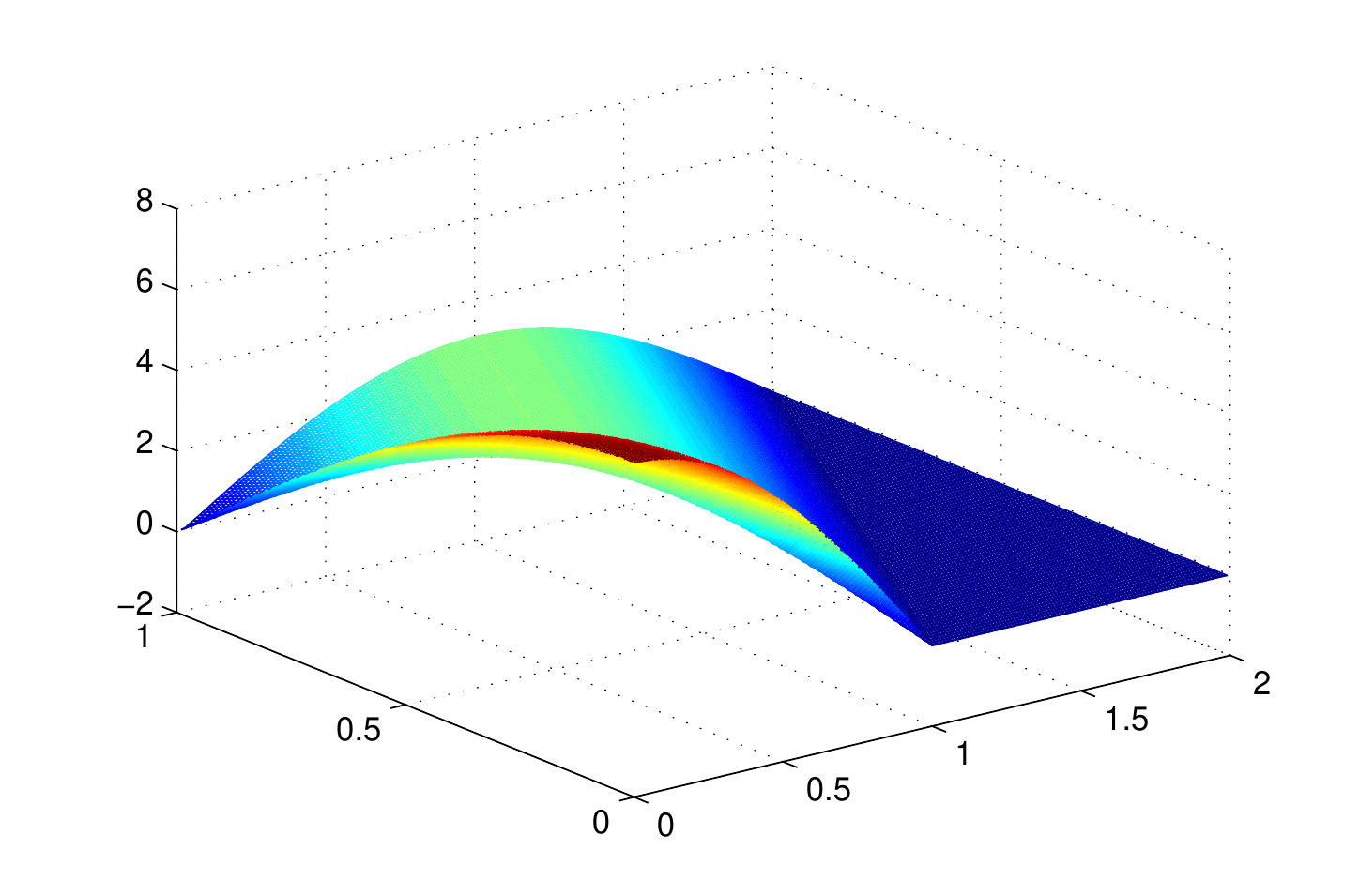}

State $y(t,x)$ and $y_x(t,x)$ with optimal {\sc Neumann} boundary
control, $T=2$.

The state is continuous.


} 


\frame{ \frametitle{
Example: Control time $T=10$}

Optimal state for the  control time $T=10$:
%
\includegraphics[width=5cm]{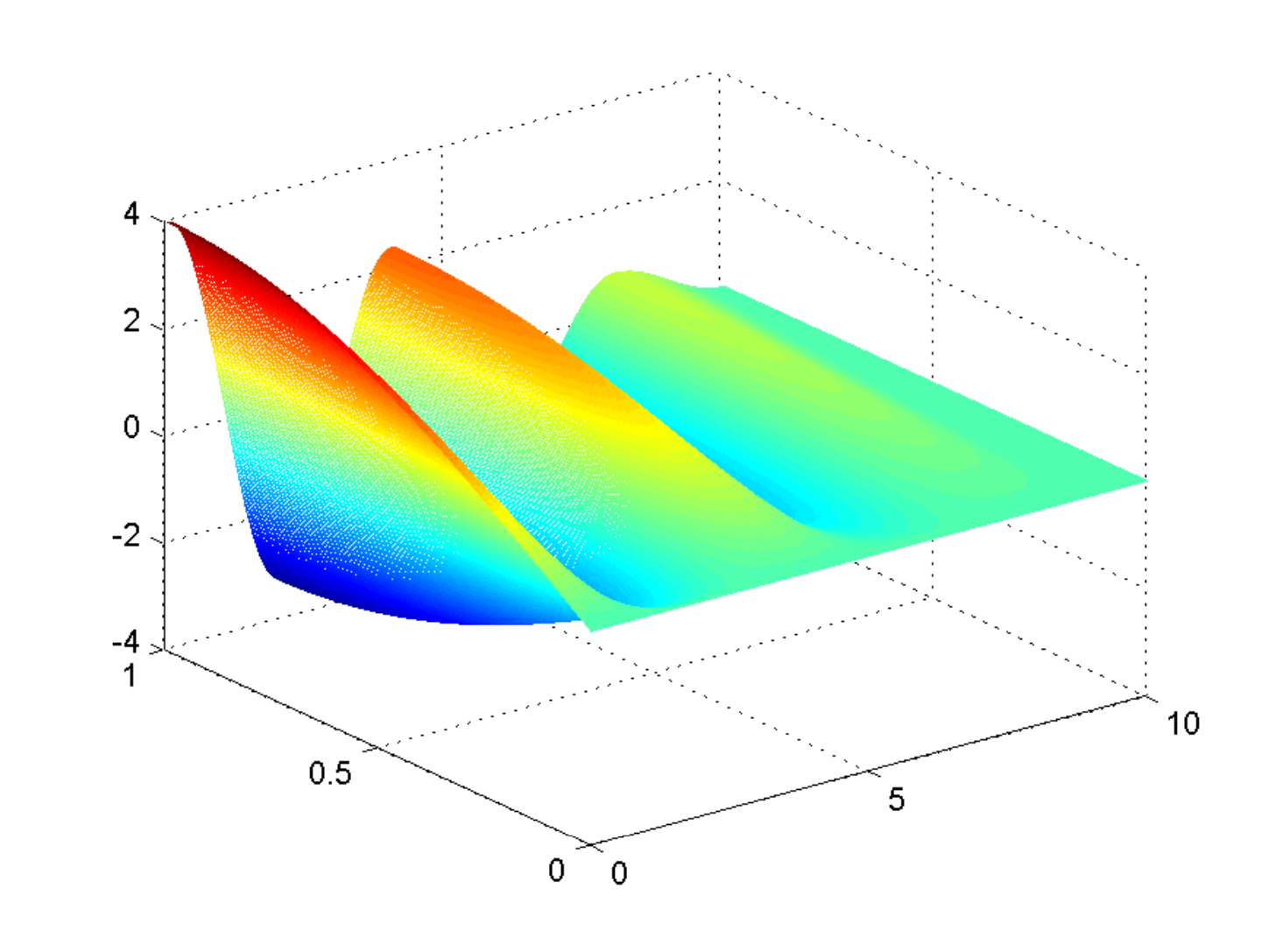}
\includegraphics[width=5cm]{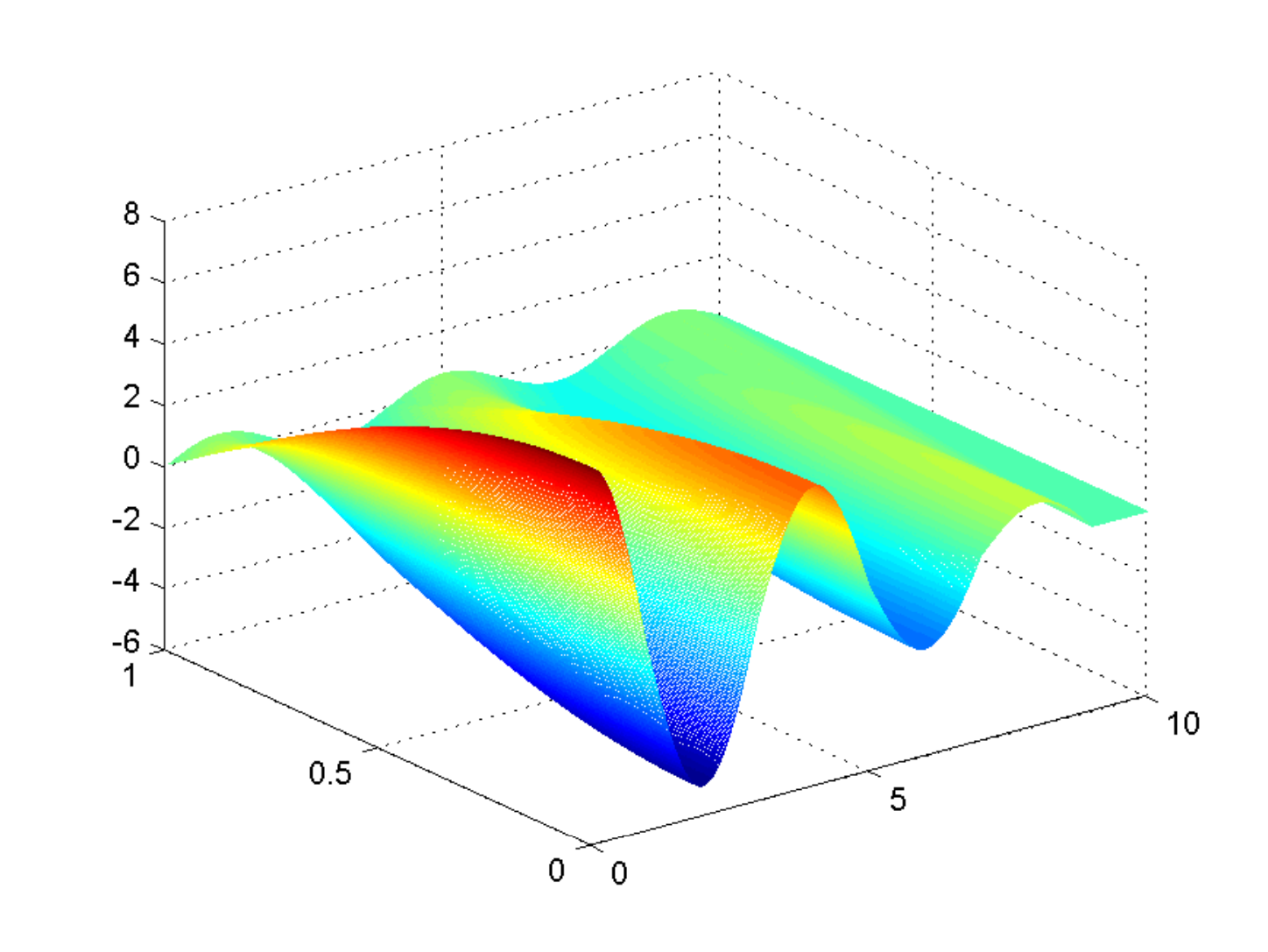}

State $y(t,x)$ and $y_x(t,x)$ with optimal {\sc Neumann} boundary
control, $T=10$.



} 

\frame{ \frametitle{
Example: Control time $T=20$}

Optimal state for the  control time $T=20$:
%
\includegraphics[width=5cm]{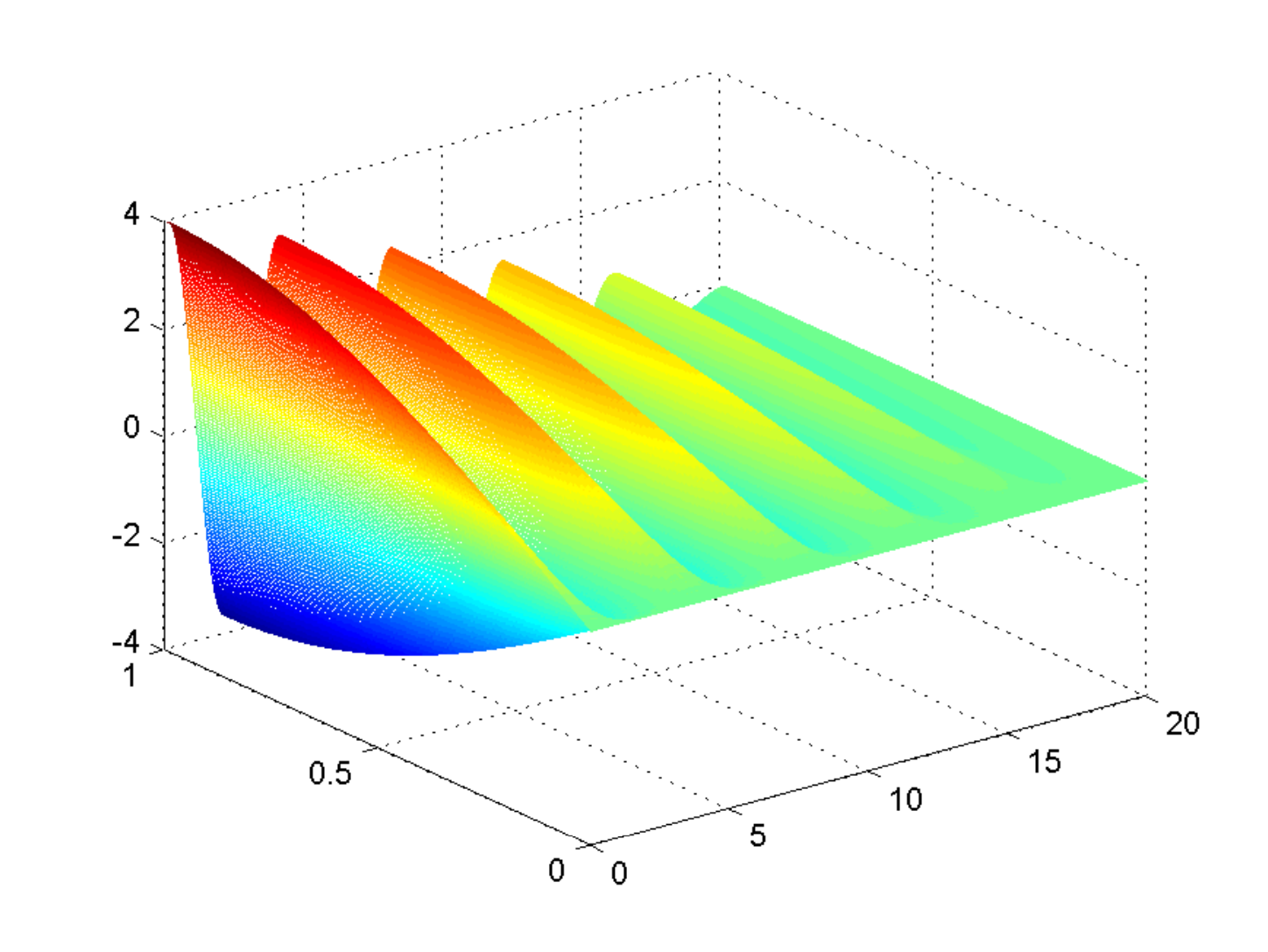}
\includegraphics[width=5cm]{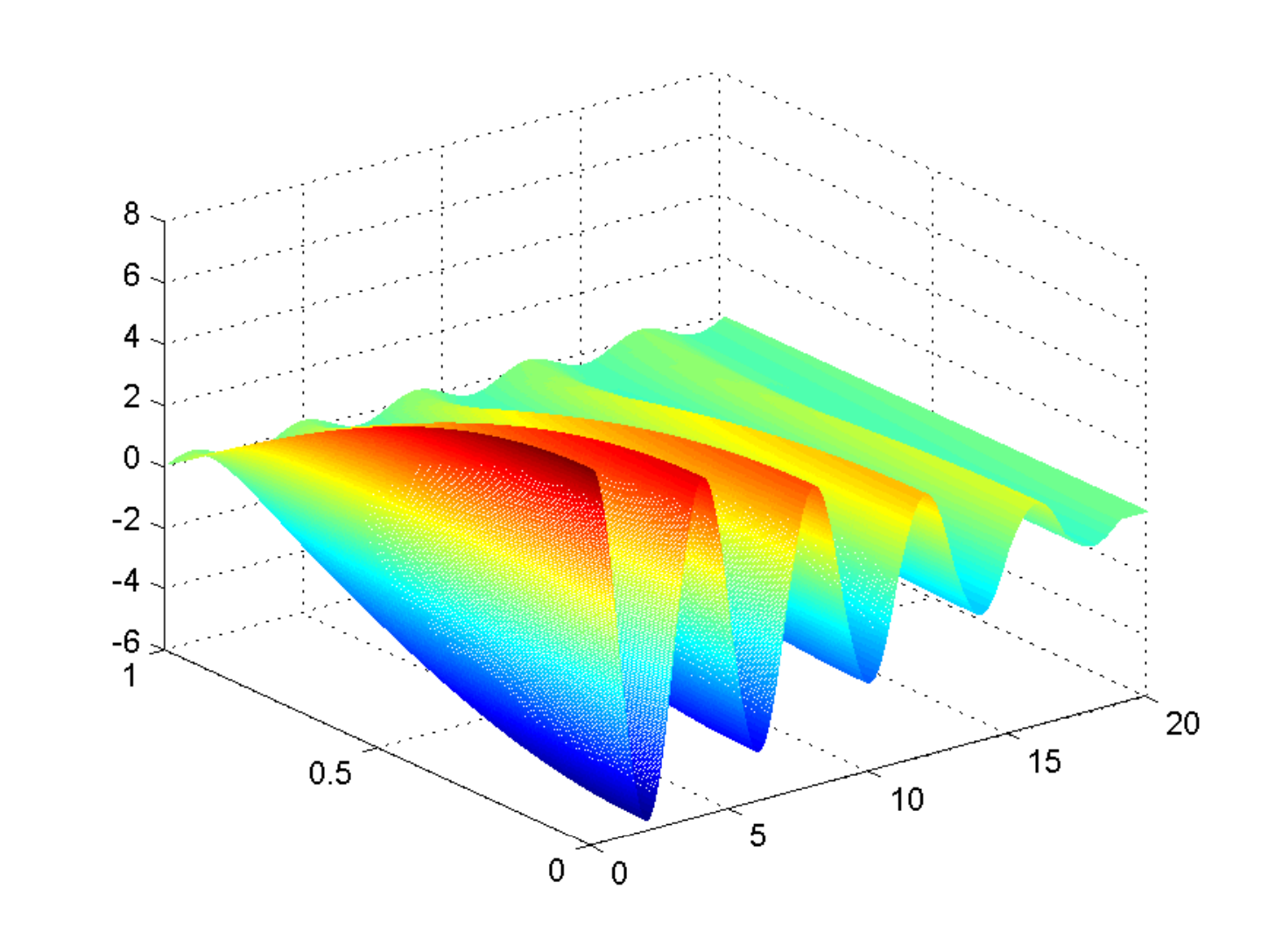}

State $y(t,x)$ and $y_x(t,x)$ with optimal {\sc Neumann} boundary
control, $T=20$.



} 

\frame{

\begin{itemize}

\item<1->
The optimal exact control sucks the energy out of the system in finite time $T$.

\item<2->

For $T=2n$, we have $n$ time subintervals of equal length.
\\
 In each subinterval
the $n$th part of the initial energy is taken out of the system.

\item<3->
Now we  look at 
{\huge stabilization} where in general, we never reach zero energy.

\end{itemize}

} 















\section{Stabilization}


\frame{ \frametitle{Stabilization}

\begin{itemize}
\item<1->

The open loop control depends on the initial state $(y_0,\,y_1)$.

In general, this state is {\bf not} known.



What happens, if the true initial state is a different from $(y_0,\,y_1)$?

\item<2->
Example: $\tilde y_0(x)= 2x$, $y_1(x)=0$.

\includegraphics[width=9cm]{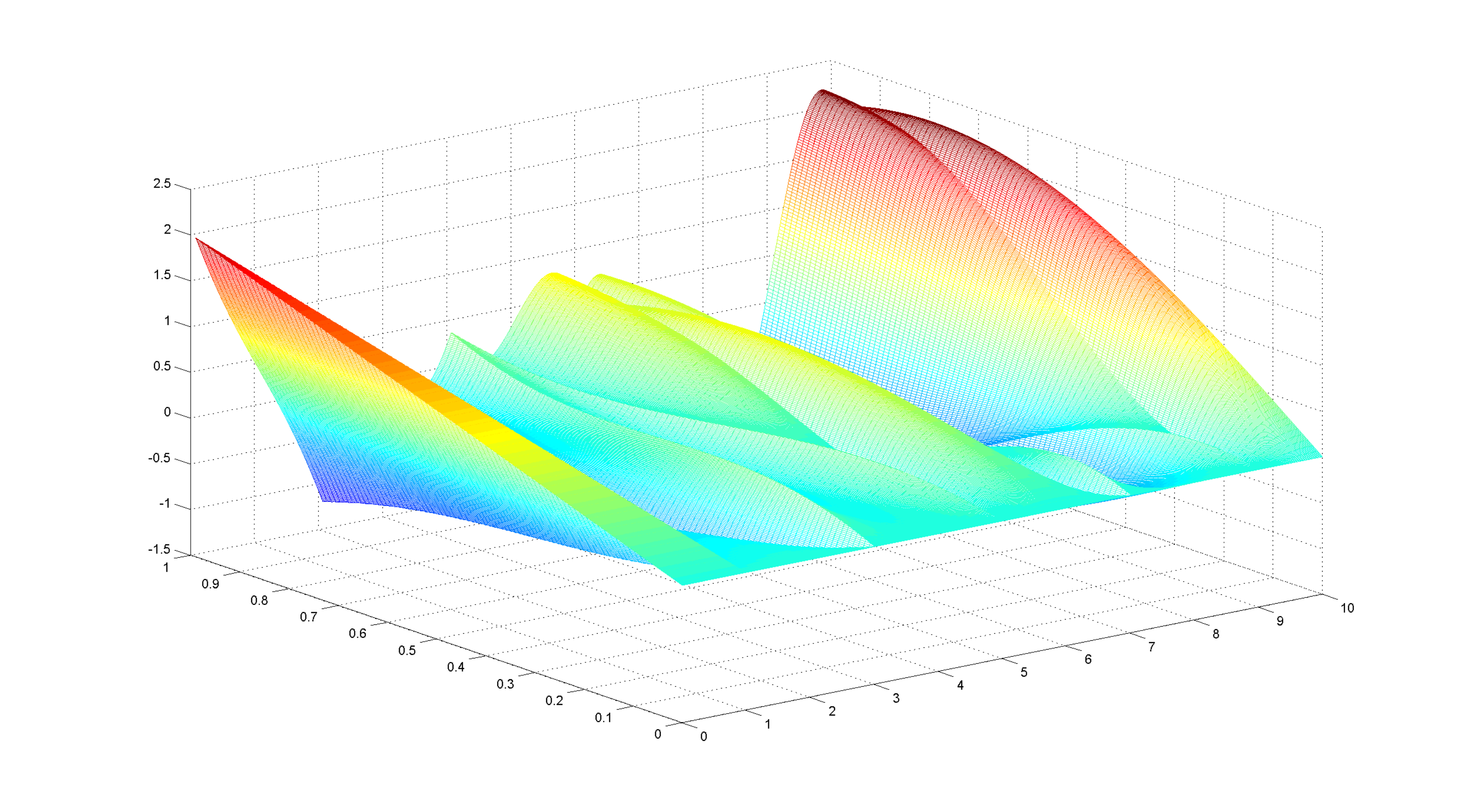}

\end{itemize}

} 


\frame{ \frametitle{Stabilization}

\begin{itemize}
\item<1->
Let $f$ be a real number. This is our feedback parameter.

\item<2->


Introduce a feedback law (closed loop control) at $x=1$:


 $$y_x(t,1)= - f \; y_t(t,1)$$

\item<3->

$$
\label{STAB} {\bf (STAB)}
\left\{
\begin{array}{l}
y(0,x)= \tilde y_0(x),\; y_t(0,x)= \tilde y_1(x),\; x\in (0,1)
\\
\\
y(t,0)= 0,\; \fbox{$y_x(t,1)= - f \, y_t(t,1)$},\; t\in (0,T)
\\
\\
y_{tt}(t,x)= y_{xx}(t,x),\;(t,x)\in (0,T) \times (0,1)
\end{array}
\right.
$$

\end{itemize}

}  




\frame{ \frametitle{Exponential Stability of the System}

\begin{itemize}

\item

We consider the {\bf Energy}
$$E(t)= \frac{1}{2} \int_0^1
\left(y_x(t,x)\right)^2 +  \left( y_t(t,x) \right)^2 \,dx.
$$

\item

For all $f>0$
System ${\bf STAB}$ is {\em exponentially stable}, that is
there exist
$C_1,\,\mu \in (0,\infty)$
such that
$$ E(t) \leq C_1 \;
E(0)\,\exp(-\mu t), \; (t \in [0,\,\infty)).
$$



\item
 For $f=1$ ${\bf STAB}$ satisfies $y(2,x)=y_t(2,x)=0$,
\\
for all initial states! ({\em Komornik, Cox and Zuazua})

\end{itemize}

}  


\subsection{Example: Stationary Feedback Law}



\frame{ \frametitle{Example: Feedback}

\begin{itemize}

\item<1->Feedback switched off ${\bf f=0}$ (Conservation of energy):

\includegraphics[width=4cm]{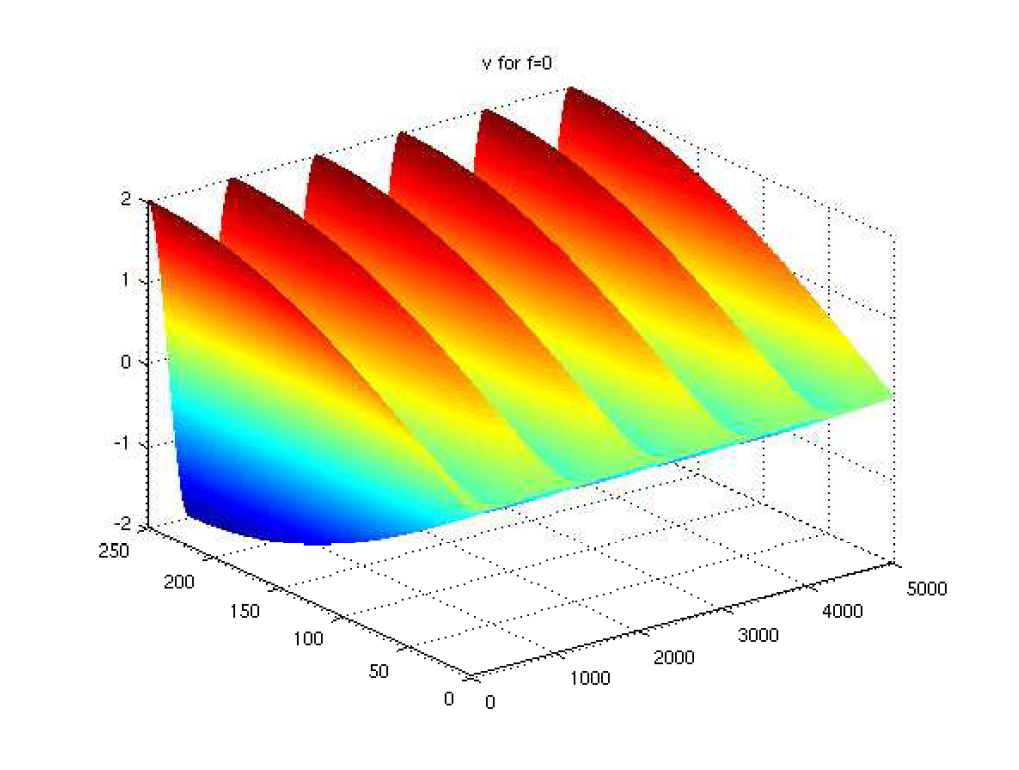}

$y(t,x)$ with $f=0$, Zero control

\item<2->
Feedback with ${\bf f=1}$:

\includegraphics[width=4cm]{bild2T=2}

State $y(t,x)$
with
feedback for $y_0=4 \sin(\pi x/2)$, $y_1=0$


\end{itemize}

} 

\subsection{Example: Time-varying Feedback Control}

\frame{\frametitle{Example: Combination $y_x = -y_t + u$}

\begin{itemize}
\item<1->
{\bf Example}
State for the control time $T=10$
with $f=1$ and the optimal control from
{(\bf EC}) for $y_0=4\sin(\frac{\pi}{2} x)$, $y_1(x)=0$
with $\tilde y_0(x)= 2x$, $y_1(x)=0$.

$$y_x(t,1) = -y_t(t,1) + u(t)$$

%
\includegraphics[width=7cm]{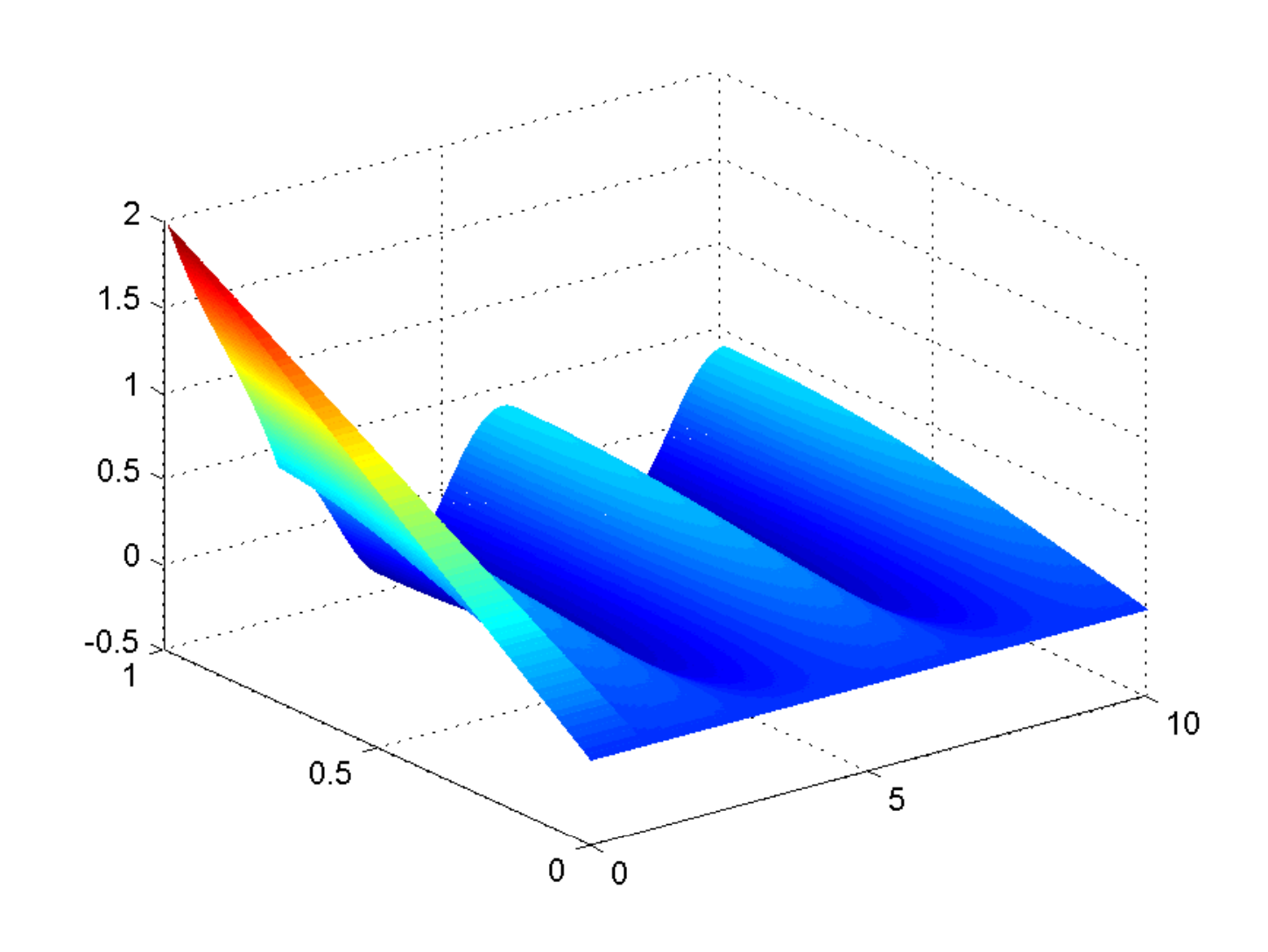}

state $y(t,x)$
with  {\sc
Neumann}-boundary control $y_x=-y_t + u_0$, $T=10$

\item<2->
Can we do better?


\end{itemize}

} 

\section{Optimized Feedback}




\frame{ \frametitle{Optimized  Feedback}

\begin{itemize}
\item<1->
To guarantee stability of the system also if
an optimal control is used,
we look at {\bf optimized Feedback}.

Let a feedback parameter $f\geq 0$ be given.


\item<2->

$$
\label{OF} {\bf (OF)}
\left\{
\begin{array}
{l}
{\rm minimize}_{u\in L^2(0,T)} \;
 \|y_x(t,1)\|_{L^2(0,T)}^2
\; {\rm subject } \;{\rm to}\;
\\
\\
y(0,x)= y_0(x),\; y_t(0,x)= y_1(x),\; x\in (0,1)
\\
\\
y(t,0)= 0,\; \fbox{$y_x(t,1)= - f y_t(t,1) + u(t)$,}\; t\in (0,T)
\\
\\
y_{tt}(t,x)= y_{xx}(t,x),\;(t,x)\in (0,T) \times (0,1)
\\
\\
y(T,x)= 0,\; y_t(T,x)= 0,\; x\in (0,1).
\end{array}
\right.
$$

\item<3->
For $f=0$ we get again {\bf (EC)}.

\item<4->
Here the optimal control depends on $y_0$, $y_1$ and $f$.

\item<5->
Due to the  objective function, the optimal value is independent of
$f$.

\item<6-> After time $T$ the control $u$ is switched off:
 $u(t)=0$ for $t>T$.
 This yields exponential stability of the system.

\end{itemize}
}  

\frame{ \frametitle{The optimal control}

\begin{itemize}

\item<1->


{\bf Theorem}
[{\small \em Gugat 2013}]
Let $T=K+1$ be even.

\item<2->
Then the optimal control
%
for $k\in \{0,1,...,(K-1)/2\}$, $t\in (0,2)$
is:
$$
u(t+2k) =
\left\{
\begin{array}{ll}
\frac{(-1)^k}{T}
\left[ 1  - f \left( T - (2k+1) \right) \right] \, 2\beta'(1-t), &  t\in (0,1)
\\
\\
\frac{(-1)^k}{T}
\left[ 1 - f  \left( T - (2k+ 1) \right) \right] \, 2\alpha'(t-1), & t\in (1,2).
\end{array}
\right.
$$


\item<3->
For the minimal control time $T=2$ we get
$$
u(t) =
\left\{
\begin{array}{ll}
\left[ 1-f \right] \, \beta'(1-t), &  t\in (0,1)
\\
\\
\left[ 1  - f \right] \, \alpha'(t-1), & t\in (1,2).
\end{array}
\right.
$$

\onslide<4->

In particular for $f=1$ we get  \fbox{$u(t)=0$}.

\vspace*{5mm}

{\bf In this case the feedback law already yields the optimal control!}

\end{itemize}

} 

\subsection{Examples for optimized feedback}


\frame{ \frametitle{Example: Minimal Control Time $T=2$}

\begin{itemize}

\item<1->
State $y$
for ${\bf f=0}$ and the optimal control from {(\bf EC}) for
$y_0=4\sin(\frac{\pi}{2} x)$, $y_1(x)=0$ with $\tilde y_0(x)= 2x$,
$y_1(x)=0$

\begin{figure}[H]
\includegraphics[width=5cm]{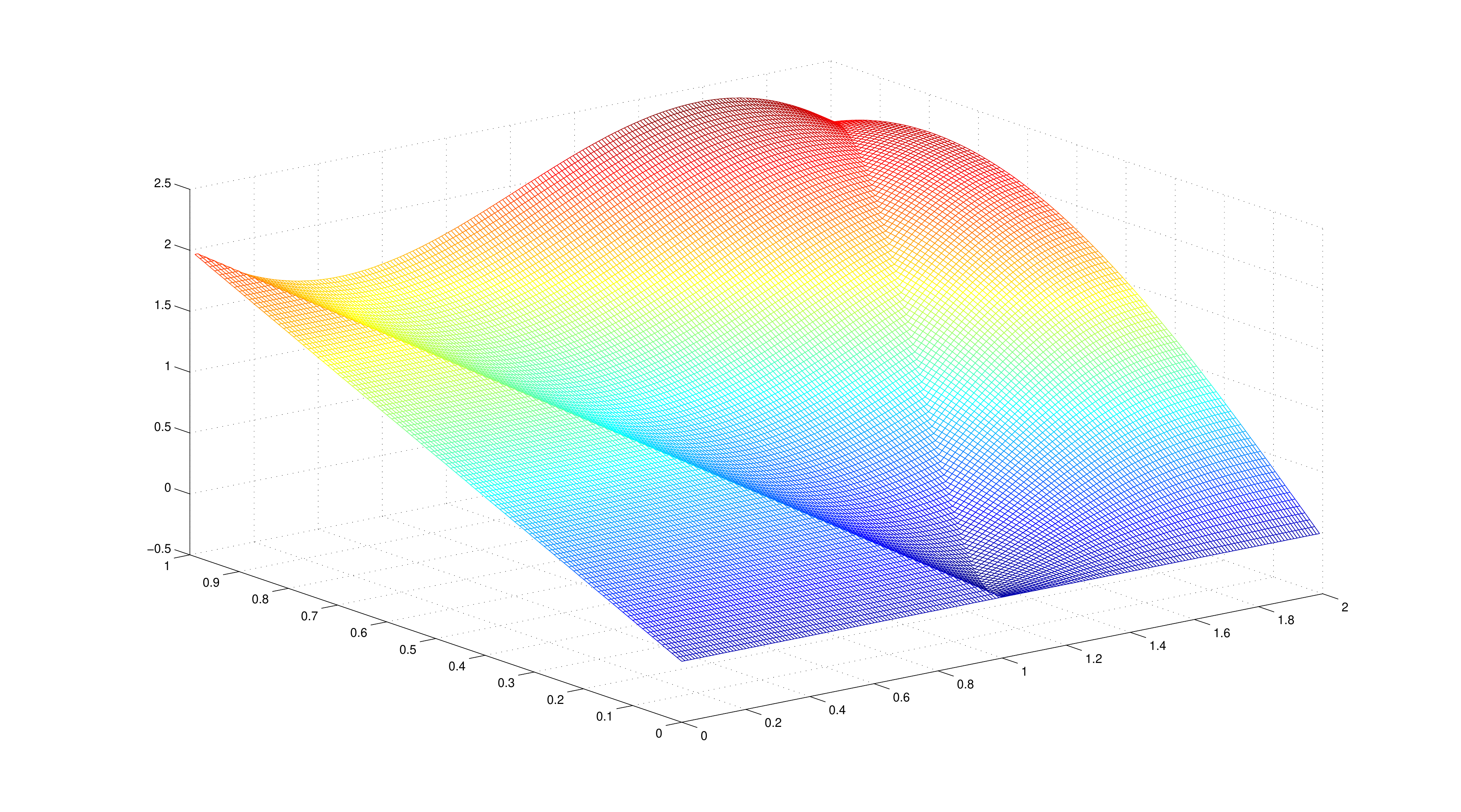}
%
\end{figure}

\item<2->


State $y$
with ${\bf f=1}$ and the optimal control $u=0$ from {(\bf OF})



\begin{figure}[H]
\includegraphics[width=5cm]{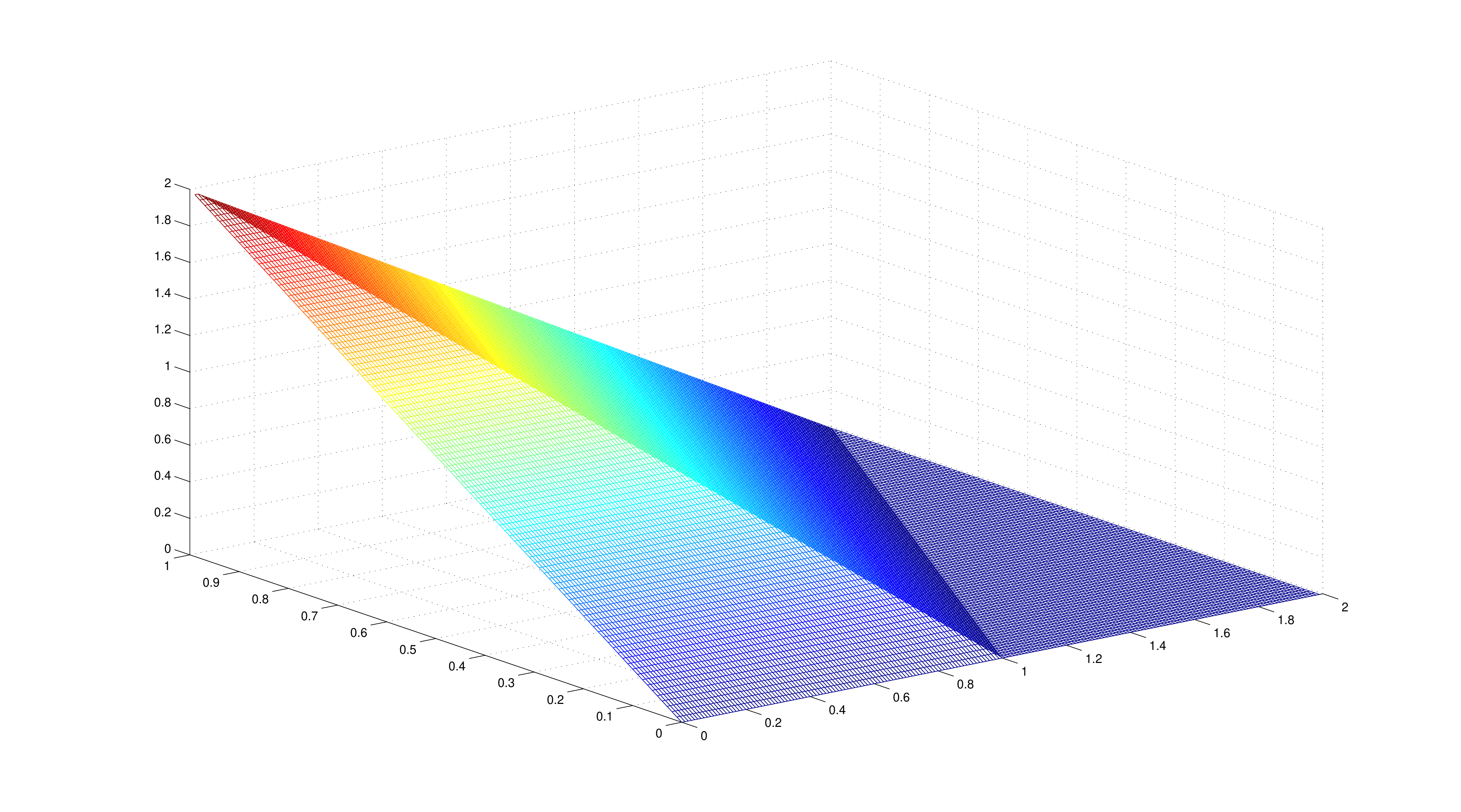}
\end{figure}

\end{itemize}

} 


\frame{ \frametitle{Example: Minimal Control Time $T=2$}

\begin{itemize}

\item



state $y$
with ${\bf f=\tfrac{1}{2} }$ and the optimal control from {(\bf
OF})



\begin{figure}[H]
\includegraphics[width=5cm]{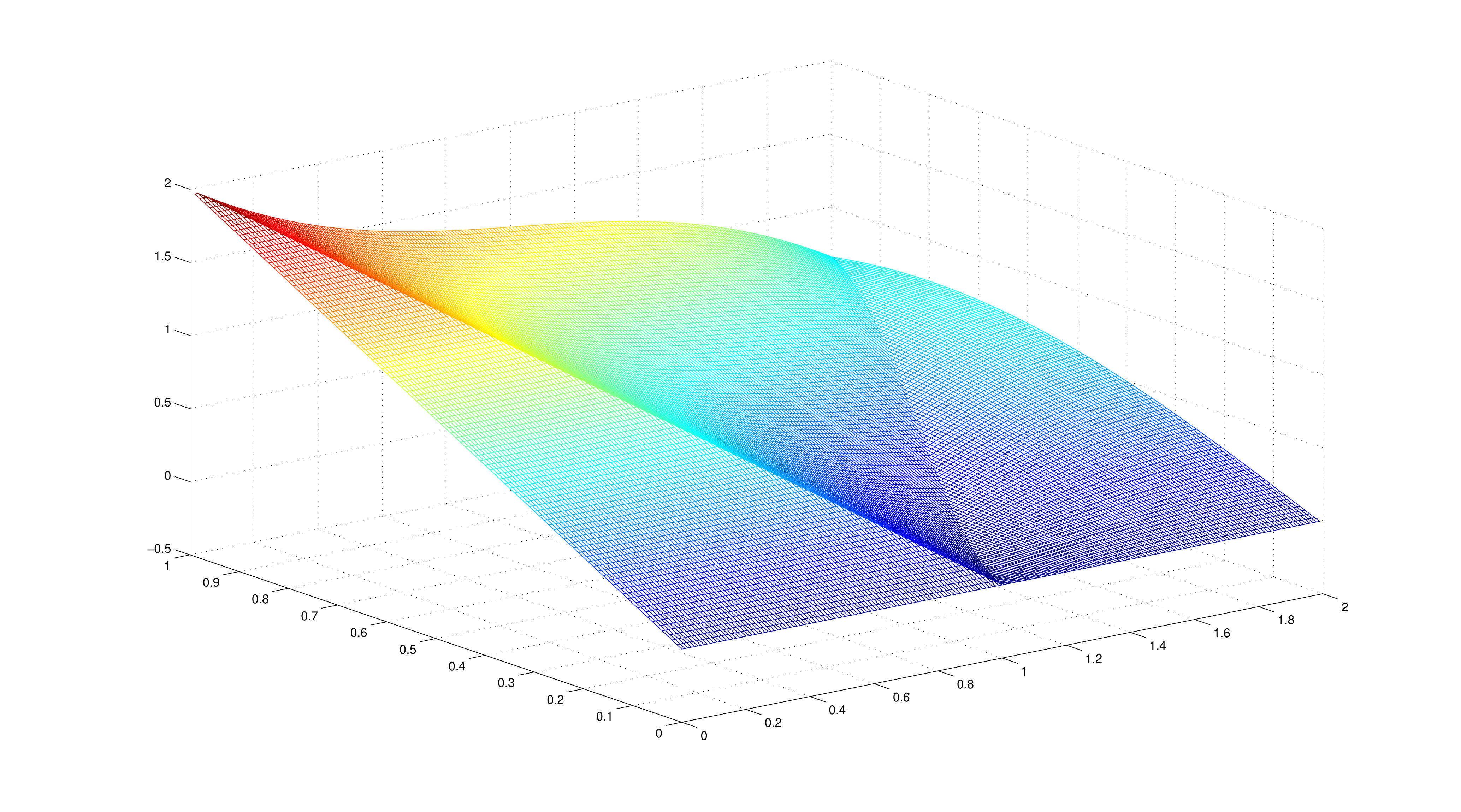}
\end{figure}

\item



State $y$
with ${\bf f=2}$  and the optimal control from
{(\bf OF})



\begin{figure}[H]
\includegraphics[width=5cm]{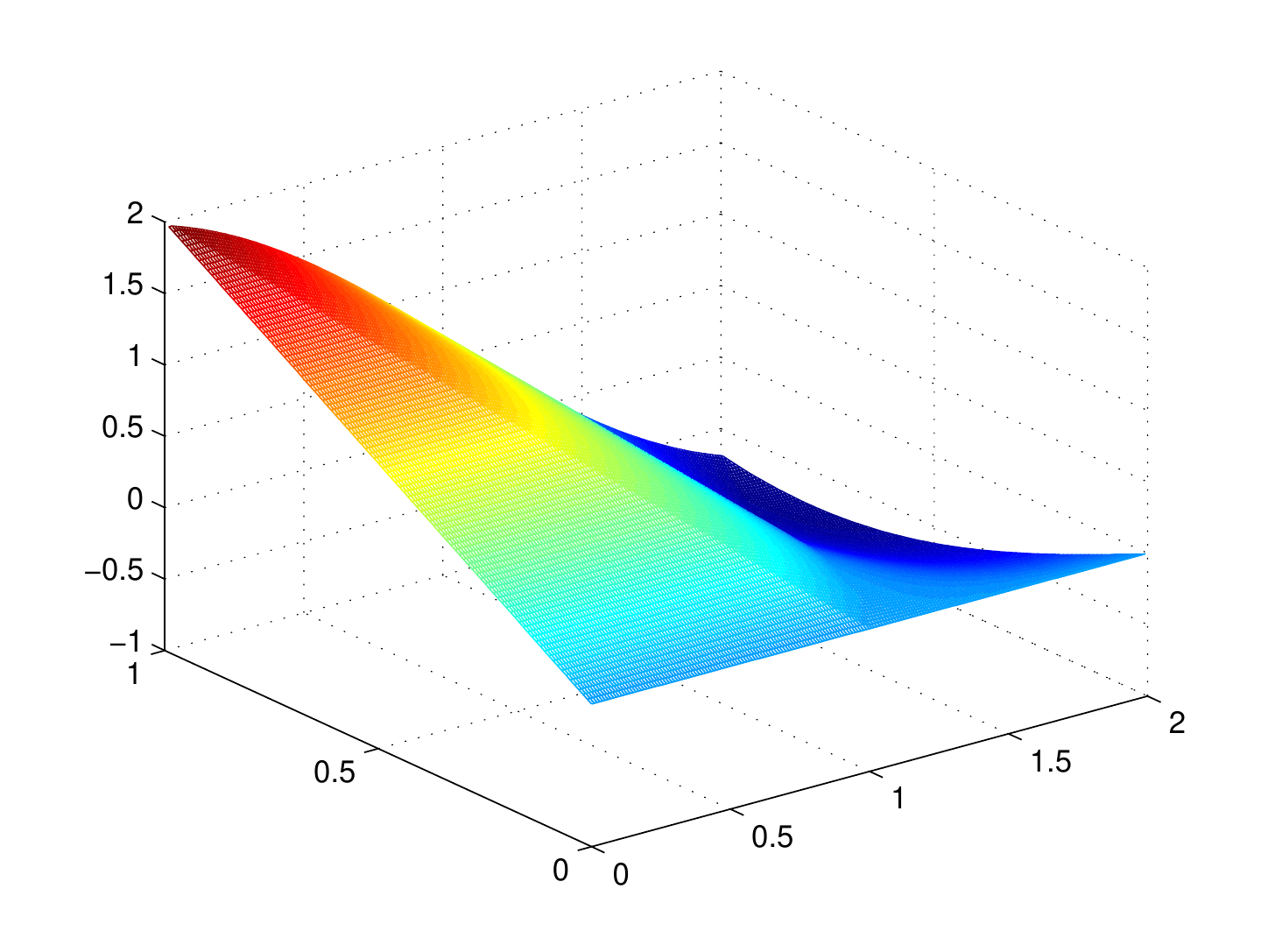}
\end{figure}

\item

With initial state
$y_0=4\sin(\frac{\pi}{2} x)$, $y_1(x)=0$
the picture is independent of $f$!
\end{itemize}

} 

\frame{ \frametitle{The optimal control in $[T-2,T]$}

\begin{itemize}

\item<1->Let  $T=K+1$ be even.

For $k\in \{0,1,...,(K-1)/2\}$, $t\in (0,2)$ we have the optimal control
$$
u(t+2k) =
\left\{
\begin{array}{ll}
(-1)^k
\frac{2}{T} \left[ 1 - f \left( T - (2k+1)\right) \right] \, \beta'(1-t), &  t\in (0,1)
\\
\\
(-1)^k
\frac{2}{T} \left[ 1 - f \left( T - (2k+1)\right) \right]
\, \alpha'(t-1), & t\in (1,2).
\end{array}
\right.
$$

\item<2->
For $2k=T-2$ this implies
$$
u(t+T-2) =
\left\{
\begin{array}{ll}
(-1)^k
\frac{2}{T} \left[ 1 - f  \right] \, \beta'(1-t), &  t\in (0,1)
\\
\\
(-1)^k
\frac{2}{T} \left[ 1 - f\right]
\, \alpha'(t-1), & t\in (1,2).
\end{array}
\right.
$$

\item<3->
%
%
%
%

Hence for $f=1$ the optimal control satisfies

%
$$\fbox{$u(t)|_{[T-2,T]}=0$}.$$
With $f=1$ and $u$  with all initial states at time  $T$
the zero state is reached {\bf exactly}!

\end{itemize}

} 

%
%
%
%
%
%
%
%
%
%
%
%
%
%
%

%
%
%
%
%
%
%
%
%
%
%
%
%
%
%
%


\frame{ \frametitle{Example: Control Time $T=20$}

\begin{itemize}

\item<1->
The optimal control from  {(\bf OF)} for $T=20$ and $f=1$

\includegraphics[width=7cm]{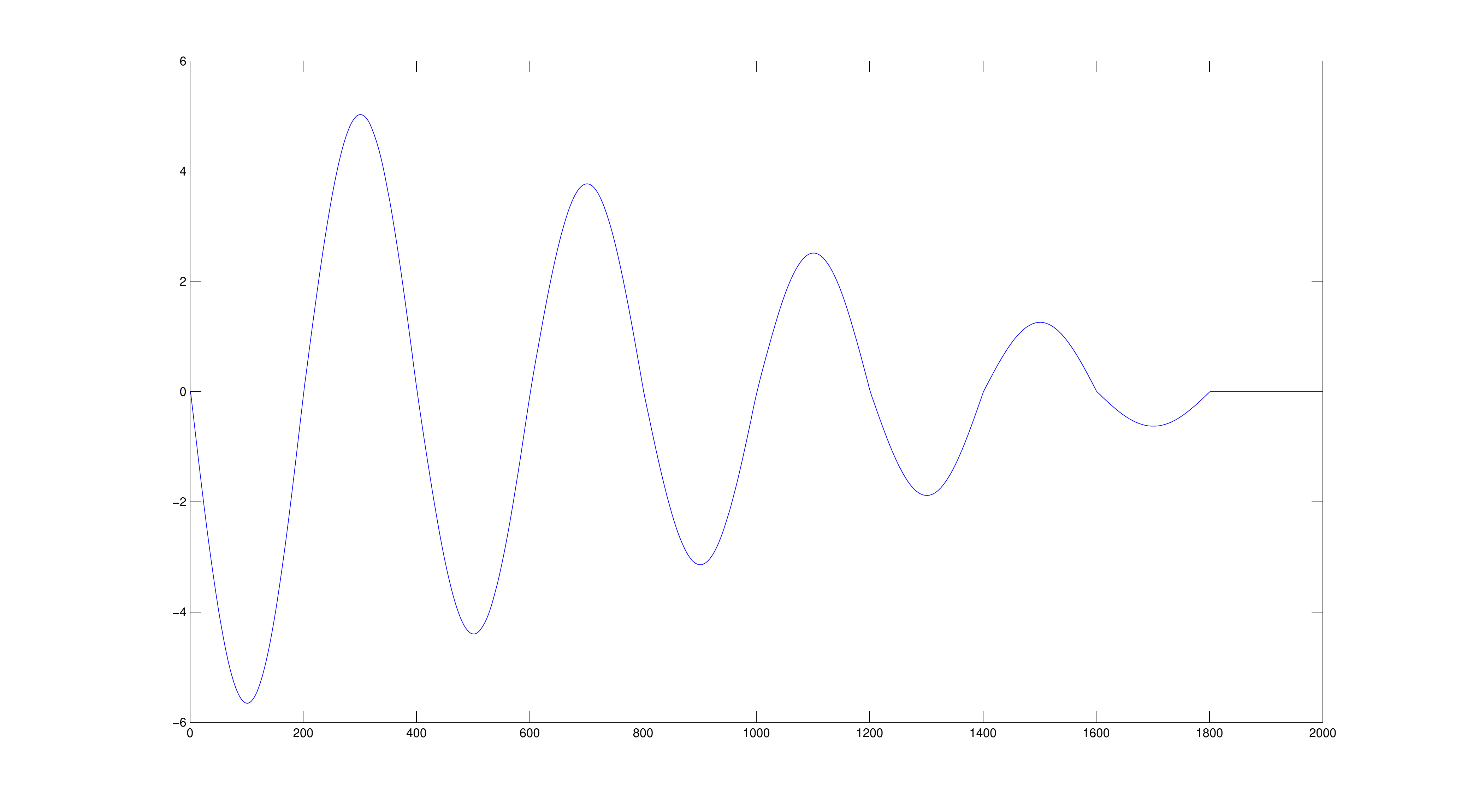}

\item<2->

The generated state with this control and intial state $\tilde y_0(x)= 2x$,
$y_1(x)=0$

\includegraphics[width=7cm]{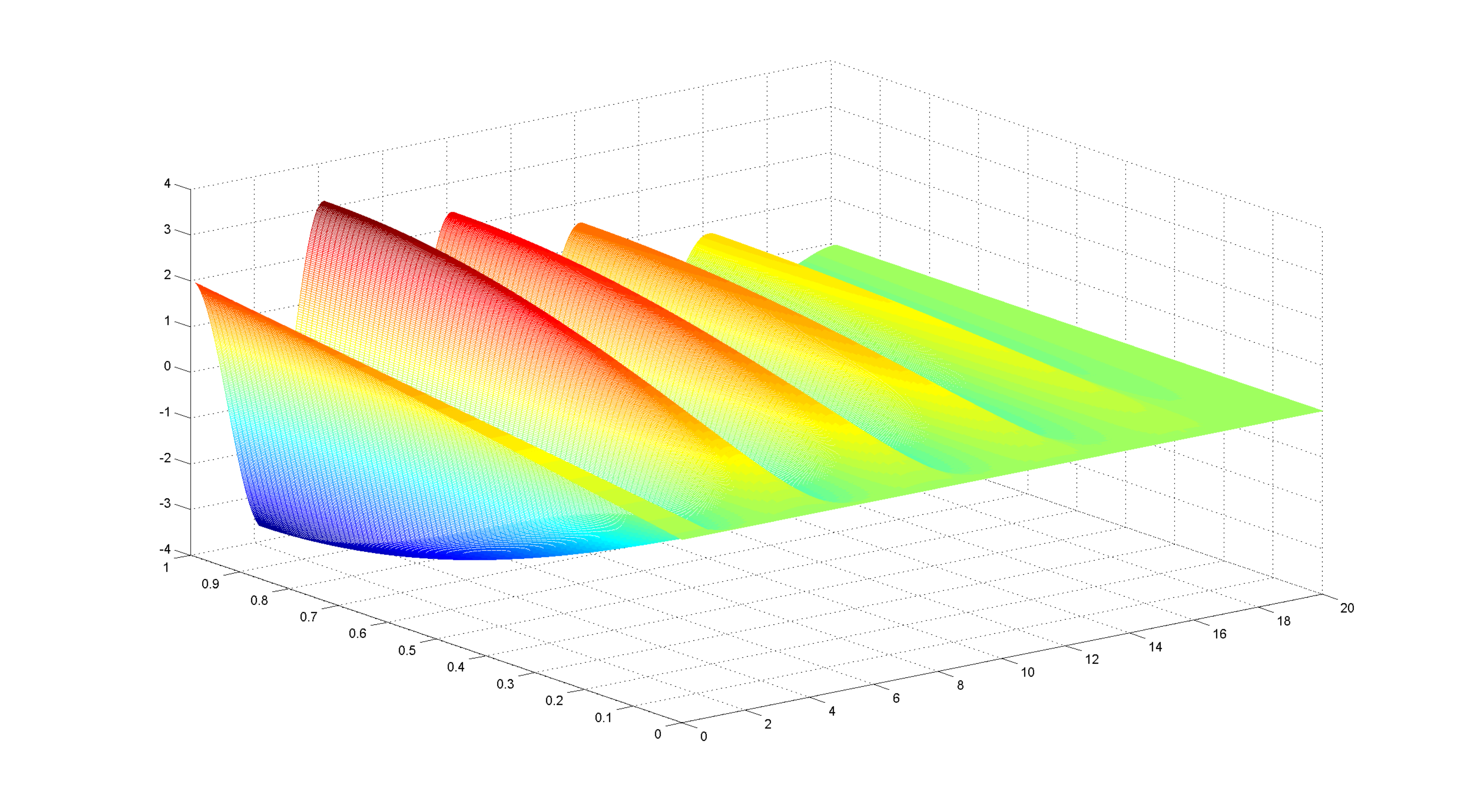}

\end{itemize}

} 


\frame{ \frametitle{Example: Control Time $T=20$}

\begin{itemize}

\item<1->
State $y$ with $\tilde y_0(x)= 2x$, $y_1(x)=0$, ${\bf f=\tfrac{1}{2}}$

\includegraphics[width=7cm]{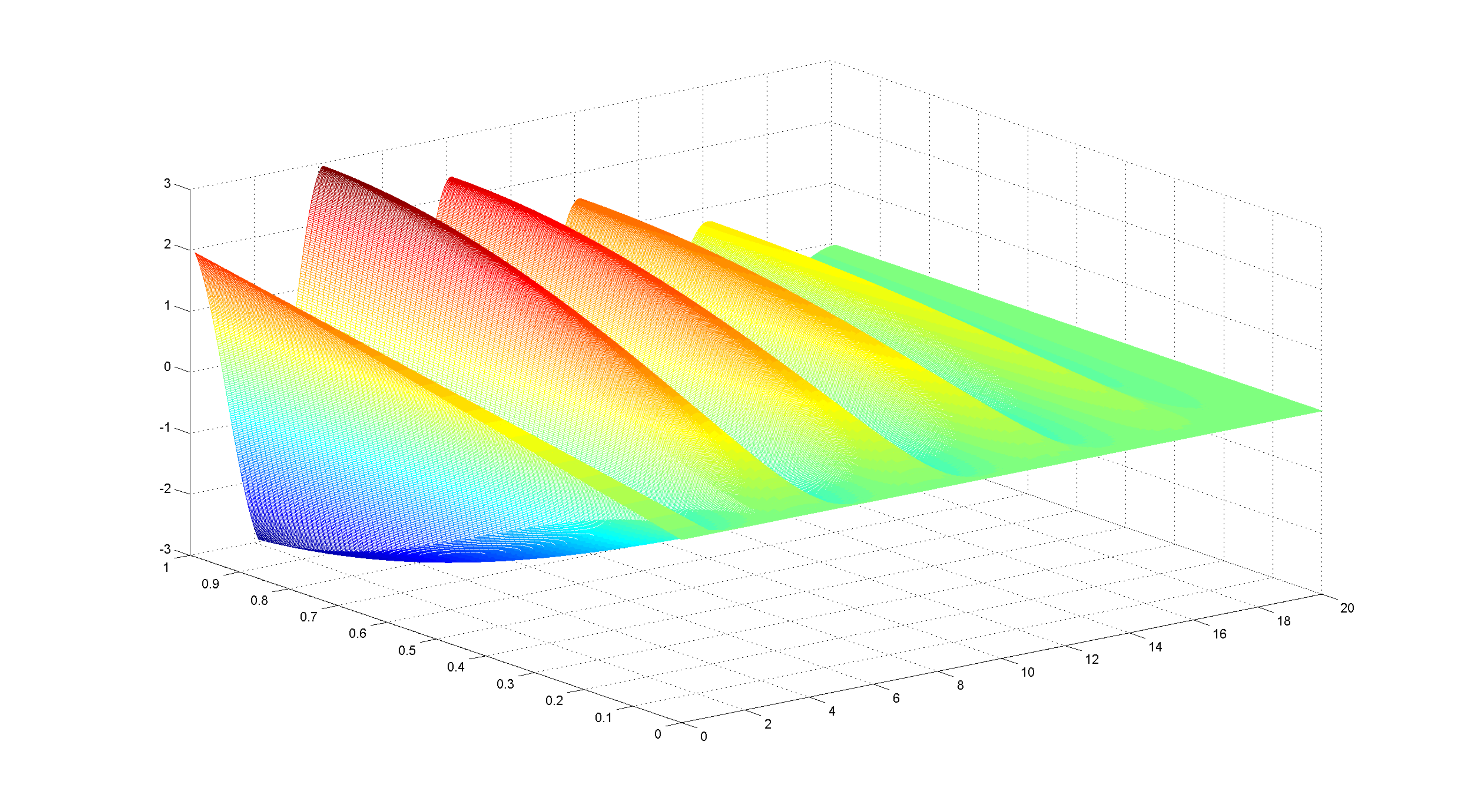}

\item<2->

State $y$ with $\tilde y_0(x)= 2x$, $y_1(x)=0$, ${\bf f=0}$  (Feedback control switched off)

\includegraphics[width=7cm]{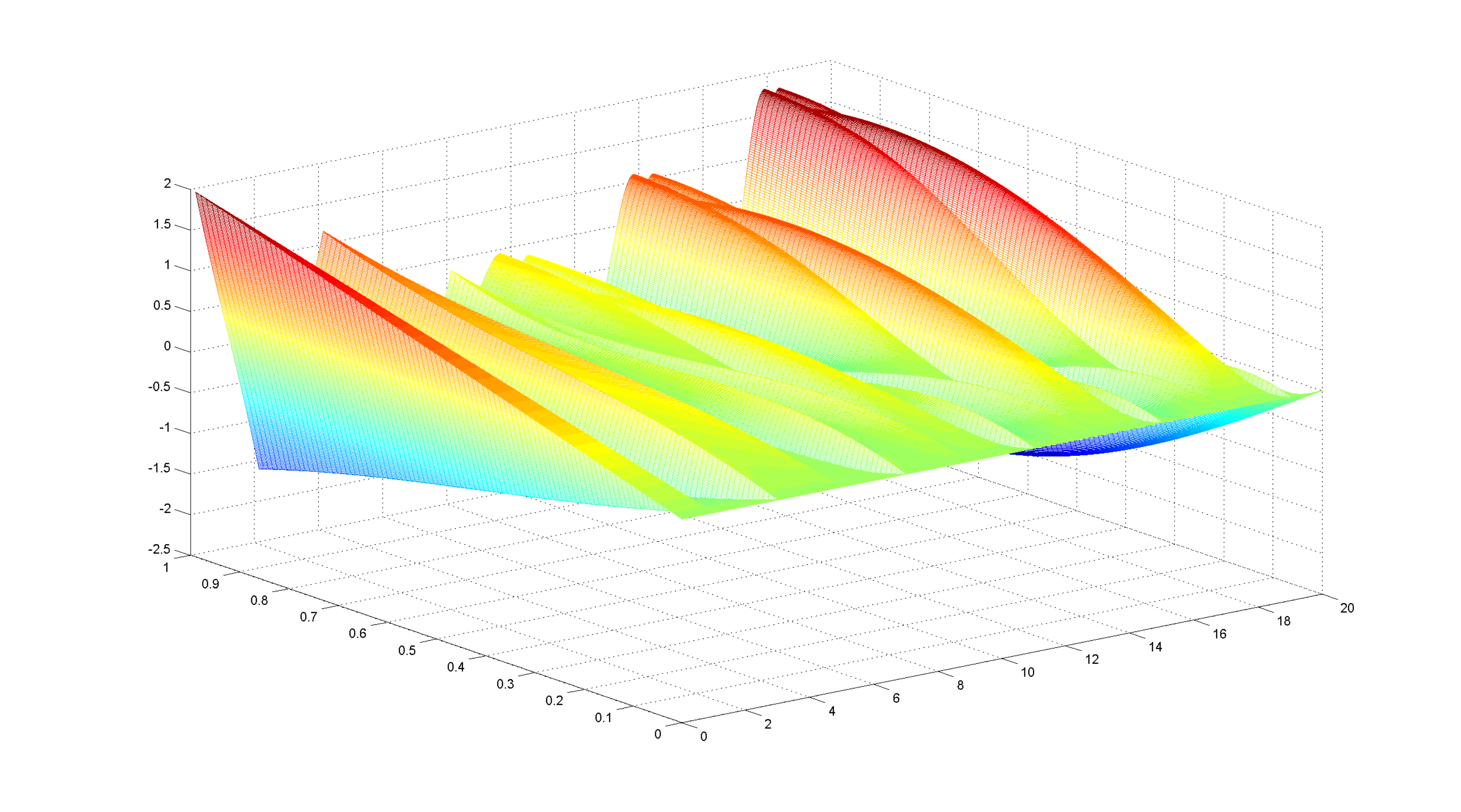}

\end{itemize}

} 



\subsection{Examples for optimized feedback:  Korteweg-de Vries Equation}


\frame{ \frametitle{Korteweg-de Vries}

{\em Cerpa and Coron 2013:}
\\
Feedback stabilization with {\bf exponential stability}
with a  suitably chosen kernel $k$
for initial state with
a sufficiently small $L^2$-norm:

$$
\left\{
\begin{array}{l}
y(0,x)=y_0(x) \in L^2(0,1)\\
y_t + y_x + y_{xxx} + y y_x =0\\
y(t,1)=0\\
y_x(t,1)=0\\
\fbox{
$y(t,0)= \int_0^1 k(0,z)  y (t,z) \, dz$}
\end{array}
\right.
$$

Method: Backstepping.
\\
For the integral feedback, the information in $y(t,z)$, $z\in (0,1)$ is used

} 


\frame{ \frametitle{Time-varying Feedback Control:  Korteweg-de Vries}

The System is locally  exactly controllable to zero.
\\
{\em L. Rosier: Control of the surface of a fluid by a wavemaker, ESAIM:COCV 10 (2004)}

{\em O.  Glass, S. Guerrero:
Some exact controllability results for the linear KdV equation and uniform controllability in the zero-dispersion limit, Asympt. Anal.  (2008)}

{\bf Optimized Feedback stabilization (with respect to $y_0$)}

$$
\left\{
\begin{array}{l}
\inf_u
\int_0^T
y(t,0)^2
\, dt
\; {\rm subject } \;{\rm to}\;\\
y(0,x)=y_0(x) \in L^2(0,1)\;  {\rm small}\\
y_t + y_x + y_{xxx} + y y_x =0\\
y(t,1)=0\\
y_x(t,1)=0\\
\fbox{
$y(t,0)= \int_0^1 k(0,z)  y (t,z) \, dz + u(t)$}
\\
y(T,x)=0.
\end{array}
\right.
$$

%
%
%
%
%

If the inital state $y_0$ is known exactly (which is never the case), this gives
exact control to zero. Otherwise exponential stability (with $u(t)=0$ for $t\geq T$).

} 



\frame{ \frametitle{Optimized Feedback Control:  Korteweg-de Vries}

\begin{itemize}

\item<1->
{\bf Step 1:}
From  ({\em Glass, Guerrero}):
Determine an exact control $v$ that is $\varepsilon$-optimal/feasible
 for

$$
\left\{
\begin{array}{l}
\inf_v
\int_0^T  (v(t))^2\, dt
\; {\rm subject } \;{\rm to}\;\\
y(0,x)=y_0(x) \in L^2(0,1)\;  {\rm small}\\
y_t + y_x + y_{xxx} + y y_x =0\\
y(t,1)=0\\
y_x(t,1)=0\\
\fbox{
$y(t,0)
= v(t)$}
\\
y(T,x)=0.
\end{array}
\right.
$$
\item<2->
{\bf Step 2:} Set
$$u(t)= v(t) - \int_0^1 k(0,z)  y_{y_0,v} (t,z) \, dz$$
 where $v(t)=0$ for $t>T$.

Then by {\em Cerpa, Coron} the system with control
\fbox{
$y(t,0)= \int_0^1 k(0,z)  y (t,z) \, dz + u(t)$}
is exponentially stable and
if $y(0,\cdot)=y_0$, it is steered to zero at time $T$.

\end{itemize}

} 



\begin{frame}

{\Huge \bf Stabilization of semilinear wave equations}

\end{frame}



\section{Stabilization of semilinear systems: Telegraph equation}


\frame{ \frametitle{Semilinear wave equation}
\begin{itemize}
\item<1-> For initial data
 $y_0\in L^{\infty}(0,1)$, $y_1\in W^{-1,\infty}(0,1)$
consider  a system with the nonlinear wave equation
(includes {\bf telegraph} equation, waterhammer eqn.)
\begin{equation}
\label{pde}
y_{tt}(t,x)  \fbox{$-  2 g_y(x, \, y (t,x)) \; y_t (t,x) $}\;  = y_{xx}(t,x)
\end{equation}
where
\begin{equation}
\label{wungleichung}
|g_y(x,\, y)| \leq w
\end{equation}
with the boundary conditions
$$y(t,0)= 0,\; \fbox{$y_x(t,1)= -  y_t(t,1)
$,}\; t\in (0,T).$$

\item<2->
For
$w<1/20$,  $\|y(t,\cdot)\|_{L^\infty(0,1)}$ decays exponentially with rate
$$\mu = \left| {\rm ln}(20 w)\right|.$$
\\
\onslide<3->
Thus the decay rate becomes arbitrarily large for $w\rightarrow 0$.

\item<4->
%
%

Consider now stability of ISS type (see {\em Mazenc, Prieur, MCRF 1, 2011}).

\end{itemize}

} 


\frame{ \frametitle{Semilinear wave equation: ISS stability}
\begin{itemize}
\item<1-> For initial data
 $y_0\in L^{\infty}(0,1)$, $y_1\in W^{-1,\infty}(0,1)$
consider  a {\bf perturbed} system
\begin{equation}
\label{pdeh}
y_{tt}(t,x) -  2 g_y(x, \, y (t,x)) \; y_t (t,x)   = y_{xx}(t,x) \;\;\fbox{$+ \;\; D(t,x)$}
\end{equation}
with continuous uniformly bounded $D$ and
$\left(
|g_y(x,\, y)| \leq w
\right)
$
with the boundary feedback
$y(t,0)= 0,\; \fbox{$y_x(t,1)= -  y_t(t,1)$}$

\item<2-> {\bf Related:} For the linear wave equation $g=0$ in
 {\em Gugat, Tucsnak,
Sigalotti: Robustness analysis for the boundary control of the
string equation}, 2007)  the influence of the position coefficient
$b$ in the feedback

 \fbox{$y_x(t,1)= - f y_t(t,1) - b y(t,1)$}

on the robustness is studied:
\\
In some cases with $b>0$, the system is more robust with respect to
$D$ than for $b=0$.

%
%
%
%
%

\end{itemize}

} 



\frame{ \frametitle{Semilinear wave equation: ISS stability
($L^\infty$)}
\begin{itemize}
\item<1->
Let $\delta$ solve the  linear closed loop system
 $\delta_{tt}= \delta_{xx} + D$, $\delta(0,x) =
\delta_t(0,x)=0$, $\delta(t,0)=0$, $\delta_x(t,1)= - \delta_t(t,1)$.

\item<2->
Due to the feedback law, the solution $\delta$ has limited memory
with respect to $D$: $\delta(t,x)$ only depends on the data
$D(s,x)|_{s\in (t-4,t)}$!

\onslide<3-> This implies in particular, that $${\rm ess}\sup_t
\|\delta(t,\cdot)\|_{L^\infty(0,1)}$$ remains bounded if $D$ is
uniformly bounded.

\item<4->
We get the {\bf robustness estimate} (for $k\in\{1,2,3,...\}$
$$
{\rm ess} \sup_{s\in [2k, 2k+2]} \|y(s,\cdot)\|_{L^\infty(0,1)}
$$
$$
\leq
(20w)^k
{\rm ess}\sup_{s\in [0, 2]}
\|y(s,\cdot)\|_{L^\infty(0,1)}
+
\frac{1 - (20w)^k}{1 - 20 w}
{\rm ess}\sup_{t\in [0,2k+2]} \|\delta(t,\cdot)\|_{L^\infty(0,1)}.
$$

\end{itemize}

} 

\begin{frame}

{\Huge \bf Stabilization of quasilinear  wave equations}

\end{frame}

\section{Stabilization of a  quasilinear wave equation}




\frame{ \frametitle{Quasilinear wave equation}
\begin{itemize}
\item<1->
In a paper with {\em Leugering, Wang, Tamasoiu}, we have studied the pde
\begin{align}\label{2.21}
\tilde u_{tt}+2\tilde u\tilde u_{tx}-(a^2-\tilde u^2)\tilde
u_{xx}=\tilde F(\tilde u,\tilde u_x,\tilde u_t).
\end{align}
with {\sc Neumann} boundary control.
\item<2->
To stabilize the system governed by
the quasilinear wave equation (\ref{2.21}) locally around a
stationary state $\bar u(x)$, we use boundary feedback given by
\begin{align*}
&x=0:\tilde u_x=\bar u_x(0)+k\tilde u_t,\\
&x=L:\tilde u=\bar u(L),
\end{align*}
with a feedback parameter $k\in (0,\infty)$.
\item<3->
If $L$ is small enough,
for suitably chosen $k>0$, sufficiently small  $C^2$ solutions
$u = \tilde u - \bar u$
of the system decay
exponentially:
$$
\|(u(t,\cdot),u_t(t,\cdot))\|_{H^2(0,L)\times
H^1(0,L)}
\leq
\eta_1
\|(u(0,\cdot),u_t(0,\cdot))\|_{H^2(0,L)\times
H^1(0,L)}
\exp\left(-\bar \mu t\right)
$$

\end{itemize}

} 


\frame{ \frametitle{Quasilinear wave equation}
\begin{itemize}
\item<1->

The analysis is based upon the Lyapunov function:
{
$$
E(t)=  \int_0^L
h_1(x)
\left[
 \Big((a^2-{\tilde u}^2)u_x^2+u_t^2\Big)
 +
 \Big((a^2-{\tilde u}^2)u_{xx}^2+u_{tx}^2\Big)
 \right]
 $$
 $$
 -2 h_2(x)
 \left[
\Big({\tilde u}\,u_x^2+u_tu_x\Big)\,
+
\Big({\tilde u}\,u_{xx}^2+u_{tx}u_{xx}\Big)
\right]
\,dx
$$
}
with the exponential weights
$h_1(x) = k e^{-\mu_1 x}$, $h_2(x)=e^{-\mu_2 x}$.

\item<2->
If $ \max_{(t,x)} |u(t,x)|$ is sufficiently small,
the numbers $k$, $\mu_1$, $\mu_2$ can be chosen such that
$$\|u_x\|_{H^1(0,L)}^2 + \|u_t\|_{H^1(0,L)}^2 \leq C_0 \; E(t).$$

\end{itemize}

} 



\section{Conclusion}



\frame{ \frametitle{Conclusion}

\begin{itemize}

\item<1->
Problems of optimal exact control provide optimal controls that
should be combined with a feedback law to enhance stability.
%
%
%
%
%
%
%

\item<2->
In engineering practice, we often have nonlinear dynamics on networks:

\begin{figure}[!h]

%

There are lots of open questions!


\end{figure}

\end{itemize}

} 


\begin{frame}[t]{
Thank you for your attention!
}

{\tiny
\begin{itemize}
\item
{\textit M. Gugat, G. Leugering, G. Sklyar}, $L^p$-optimal boundary
control for the wave equation, SICON 2005 \vspace{0.1cm}
\item
{\textit M. Gugat},
Optimal boundary control of a string to rest in finite time with continuous state
ZAMM, 2006 \vspace{0.1cm}
\item
{\textit M. Gugat, G. Leugering},
$L^\infty$ Norm Minimal Control of the wave equation: On the weakness of the bang--bang principle,
ESAIM: COCV 14, 254-283, 2008
 \vspace{0.1cm}
\item
{\textit M. Gugat}, Penalty Techniques for State Constrained Optimal
Control Problems with the Wave Equation, SICON 2009 \vspace{0.1cm}
\item
{\textit M. Gugat}, Boundary feedback stabilization by time delay
for one-dimensional wave equations, { IMA Journal of Mathematical
Control and Information 2010} \vspace{0.1cm}
\item
{\textit
 M. Gugat, M. Tucsnak},
  An example for the switching delay feedback
stabilization of an infinite dimensional system: The boundary
stabilization of a string, Syst.  Cont. Let. 60, 226-230,
2011
 \vspace{0.1cm}
\item
 \textit{M. Gugat, M. Herty, V. Schleper},
    Flow control in gas networks: Exact controllability to a given demand,
Mathematical Methods in the Applied Sciences 34, 745-757, 2011
 \vspace{0.1cm}
    \item {\textit Dick, M., Gugat, M. and Leugering, G,}
A strict $H^1$-Lyapunov function and feedack stabilization for the
isothermal Euler equations with friction, \textit{Numerical Algebra,
Control and optimization}, 2011 \vspace{0.1cm}
    \item {\textit Gugat, M., Dick, M. and Leugering, G.},
Gas flow in fan-shaped networks: classical solutions and feedback
stabilization, \textit{SICON}, 2011 \vspace{0.1cm}
    \item  {\textit Gugat, M. and Herty, M.},
Existence of classical solutions and feedback stabilization for the
flow in gas networks, \textit{ESAIM  COCV}, 2011
\vspace{0.1cm}
    \item {\textit  Gugat, M., Leugering, G., Tamasoiu, S. and Wang, K.},
$H^2$-stabilization of the Isothermal Euler equations with friction:
a Lyapunov function approach,  \textit{Chin. Ann. Math}., 2012

    \item {\textit  Gugat, M., Leugering, G., Tamasoiu, S. and Wang, K.},
Boundary feedback stabilization for second-order quasilinear
hyperbolic systems: A strict $H^2$-Lyapunov function, submitted to
\textit{MCRF}, 2013

\item {\textit  Gugat, M.  Sokolowski, J.,}
A note on the approximation of Dirichlet boundary control problems for the wave equation on curved domains,
Applicable Analysis 2013.

\end{itemize}
}

\end{frame}



\end{document}